\documentclass[10pt,draft]{article}
\usepackage{amsfonts,amsmath,amssymb,cite}
\begin{document}
\title{On the derivative of the associated Legendre function of
the first kind of integer degree with respect to its order}
\author{Rados{\l}aw Szmytkowski \\*[3ex]
Atomic Physics Division, \\
Department of Atomic Physics and Luminescence, \\
Faculty of Applied Physics and Mathematics, \\
Gda{\'n}sk University of Technology, \\
Narutowicza 11/12, PL 80--952 Gda{\'n}sk, Poland \\
email: radek@mif.pg.gda.pl}
\date{\today}
\maketitle
\begin{abstract}
The derivative of the associated Legendre function of the first kind
of integer degree with respect to its order, $\partial
P_{n}^{\mu}(z)/\partial\mu$, is studied. After deriving and
investigating general formulas for $\mu$ arbitrary complex, a
detailed discussion of $[\partial
P_{n}^{\mu}(z)/\partial\mu]_{\mu=\pm m}$, where $m$ is a non-negative
integer, is carried out. The results are applied to obtain several
explicit expressions for the associated Legendre function of the
second kind of integer degree and order, $Q_{n}^{\pm m}(z)$. In
particular, we arrive at formulas which generalize to the case of
$Q_{n}^{\pm m}(z)$ ($0\leqslant m\leqslant n$) the well-known
Christoffel's representation of the Legendre function of the second
kind, $Q_{n}(z)$. The derivatives $[\partial^{2}
P_{n}^{\mu}(z)/\partial\mu^{2}]_{\mu=m}$, $[\partial
Q_{n}^{\mu}(z)/\partial\mu]_{\mu=m}$ and $[\partial
Q_{-n-1}^{\mu}(z)/\partial\mu]_{\mu=m}$, all with $m>n$, are also
evaluated.
\end{abstract}
%
%
\section{Introduction}
\label{I}
\setcounter{equation}{0}
It is the purpose of the present paper to contribute to the theory of
special functions of mathematical physics and chemistry.
Specifically, we shall add to the knowledge about the associated
Legendre functions of the first, $P_{\nu}^{\mu}(z)$, and second,
$Q_{\nu}^{\mu}(z)$, kinds (cf, e.g.,
\cite{Todh75,Ferr77,Neum78,Hein78,Hein81,Olbr87,Byer93,Hobs96,Wang04,
Barn07,Wang21,Hobs31,Snow52,Erde53,Mors53,Lens54,Robi57,Robi58,Robi59,
Jahn60,Krat60,Steg65,Magn66,MacR67,Grad94,Temm96,Prud83,Prud03,
Bryc06}). We shall touch two particular problems. First, we shall
investigate the derivative of the associated Legendre function of the
first kind of integer degree with respect to its order. Second, it
will be shown that the results of that investigation may be used to
construct, in a straightforward and unified manner, several known
representations of the associated Legendre functions of the second
kind of integer degree and order, $Q_{n}^{m}(z)$; in the past those
representations were obtained by other authors with the use of a
variety of, usually more complicated, techniques. In addition, some
possibly new expressions for $Q_{n}^{m}(z)$ will be also presented.

The literature concerning the derivative $\partial
P_{\nu}^{\mu}(z)/\partial\mu$ is very limited. Surprisingly, no
relevant expressions have been given in \cite{Bryc06}, which
otherwise contains a large collection of parameter derivatives of
various special functions. Several variants of the formula
\begin{equation}
\frac{\partial P_{\nu}^{\mu}(z)}{\partial\mu}\bigg|_{\mu=0}
=\psi(\nu+1)P_{\nu}(z)+Q_{\nu}(z),
\label{1.1}
\end{equation}
where
\begin{equation}
\psi(\zeta)=\frac{1}{\Gamma(\zeta)}
\frac{\mathrm{d}\Gamma(\zeta)}{\mathrm{d}\zeta}
\label{1.2}
\end{equation}
is the digamma function \cite{Erde53,Jahn60,Magn66,Grad94,Davi65},
may be found in \cite[page 178]{Magn66}. Robin \cite[Eq.\ (333) on
page 175]{Robi58} gave the following
representation\footnote[1]{\label{FOOT1}~The notation used in the
present work differs in some respects from that adopted by Robin in
\cite{Robi57,Robi58,Robi59}. In particular, the digamma function used
by Robin was defined as
\begin{displaymath}
\psi(\zeta)=\frac{1}{\Gamma(\zeta+1)}
\frac{\mathrm{d}\Gamma(\zeta+1)}{\mathrm{d}\zeta}
\end{displaymath}
rather than as in our Eq.\ (\ref{1.2}). Also, one should be warned
that in Eq.\ (333) on page 175 in \cite{Robi58} the factor $p!$ in
the denominator of a summand is missing.} of $\partial
P_{\nu}^{\mu}(z)/\partial\mu$:
\begin{eqnarray}
\frac{\partial P_{\nu}^{\mu}(z)}{\partial\mu}
&=& \frac{1}{2}P_{\nu}^{\mu}(z)\ln\frac{z+1}{z-1}
\nonumber \\
&& +\left(\frac{z+1}{z-1}\right)^{\mu/2}
\sum_{k=0}^{\infty}\frac{\Gamma(\nu+k+1)\psi(k-\mu+1)}
{k!\Gamma(\nu-k+1)\Gamma(k-\mu+1)}\left(\frac{z-1}{2}\right)^{k}
\qquad (|z-1|<2).
\nonumber \\
\label{1.3}
\end{eqnarray}
Brown \cite{Brow95}, being apparently unaware of the Robin's finding,
rederived the expression in Eq.\ (\ref{1.3}) in the particular case
of $\nu=n\in\mathbb{N}$ and used it to find a representation of
$Q_{n}^{m}(z)$ suitable for numerical purposes. An application of
the derivative $\partial P_{\nu}^{\mu}(z)/\partial\mu$ to the
evaluation of the integral
$\displaystyle\int_{\mathcal{C}}\mathrm{d}z\:
[P_{\nu}^{\mu}(z)]^{2}/(z^{2}-1)$, where $\mathcal{C}$ is the
contour starting from $z=0$ and returning to it after a positive
circuit round the point $z=1$, was presented by Watson \cite{Wats18}.

The detailed plan of the paper is as follows. In section \ref{II}, we
shall summarize these facts about the associated Legendre functions
of the first and second kinds of integer degrees which will find
applications in later parts of the work. In section \ref{III}, the
derivative of the associated Legendre function of the first kind,
$\partial P_{n}^{\mu}(z)/\partial\mu$ (with
$z\in\mathbb{C}\setminus[-1,1]$), is investigated. After general
considerations forming section \ref{III.1}, where no restrictions are
imposed on the order $\mu$ of the Legendre function, in sections
\ref{III.2} and \ref{III.3} we focus on the derivatives $[\partial
P_{n}^{\mu}(z)/\partial\mu]_{\mu=\pm m}$ with $m\in\mathbb{N}$.
The cases of $0\leqslant m\leqslant n$ and $m>n$ are treated
separately. A brief discussion of the derivative $\partial
P_{n}^{\mu}(x)/\partial\mu$ (with $x\in[-1,1]$) is contained in
section \ref{III.4}. Section \ref{IV} presents some applications of
the results of section \ref{III}. In section \ref{IV.1}, we obtain
several explicit representations of the associated Legendre function
of the second kind of integer degree and order, $Q_{n}^{\pm m}(z)$
and $Q_{n}^{\pm m}(x)$. Some of these representations seem to be
new, particularly these ones which generalize the well-known
Christoffel's formulas for $Q_{n}(z)$ and $Q_{n}(x)$. Finally, in
sections \ref{IV.2} and \ref{IV.3}, the derivatives
$[\partial^{2}P_{n}^{\mu}(z)/\partial\mu^{2}]_{\mu=m}$,
$[\partial Q_{n}^{\mu}(z)/\partial\mu]_{\mu=m}$, and $[\partial
Q_{-n-1}^{\mu}(z)/\partial\mu]_{\mu=m}$, all with $m>n$, are
expressed in terms of the derivatives $[\partial P_{n}^{\mu}(\pm
z)/\partial\mu]_{\mu=-m}$.

Throughout the paper, we shall be adopting the notational convention
according to which $\nu\in\mathbb{C}$, $\mu\in\mathbb{C}$,
$n\in\mathbb{N}$, $m\in\mathbb{N}$, $x\in[-1,1]$, and
$z\in\mathbb{C}\setminus[-1,1]$, the latter with
\begin{equation}
-\pi<\arg(z)<\pi,
\qquad 
-\pi<\arg(z\pm1)<\pi,
\label{1.4}
\end{equation}
hence,
\begin{equation}
-z=\mathrm{e}^{\mp\mathrm{i}\pi}z,
\qquad
-z+1=\mathrm{e}^{\mp\mathrm{i}\pi}(z-1),
\qquad
-z-1=\mathrm{e}^{\mp\mathrm{i}\pi}(z+1)
\qquad (\arg(z)\gtrless0).
\label{1.5}
\end{equation}
Furthermore, it will be understood that if the upper limit of a sum
is less by unity than the lower one, then the sum vanishes
identically.

The definitions of the associated Legendre functions of the first and
second kinds used in the paper are those of Hobson \cite{Hobs31}.
%
%
\section{The associated Legendre functions of integer degrees}
\label{II}
\setcounter{equation}{0}
The material presented in this section, to be referred to in later
parts of the paper, is based primarily on
\cite{Hobs31,Robi57,Robi58}.
\subsection{The associated Legendre function of the first kind of
integer degree}
\label{II.1}
The associated Legendre function of the first kind of a non-negative
integer degree $n$ and an arbitrary complex order $\mu$ may defined
as
\begin{equation}
P_{n}^{\mu}(z)=(\pm)^{n}\frac{1}{\Gamma(1\mp\mu)}
\frac{\Gamma(n\mp\mu+1)}{\Gamma(n-\mu+1)}
\left(\frac{z+1}{z-1}\right)^{\mu/2}
{}_{2}F_{1}\left(-n,n+1;1\mp\mu;\frac{1\mp z}{2}\right),
\label{2.1}
\end{equation}
or equivalently as
\begin{eqnarray}
P_{n}^{\mu}(z) &=& \frac{1}{\Gamma(1\mp\mu)}
\frac{\Gamma(n\mp\mu+1)}{\Gamma(n-\mu+1)}
\left(\frac{z+1}{z-1}\right)^{\mu/2}\left(\frac{z\pm1}{2}\right)^{n}
{}_{2}F_{1}\left(-n,-n\mp\mu;1\mp\mu;\frac{z\mp1}{z\pm1}\right),
\nonumber \\
\label{2.2}
\end{eqnarray}
with either upper or lower signs chosen. Both hypergeometric
functions in Eq.\ (\ref{2.1}) are proportional to the Jacobi
polynomial \cite{Magn66,Grad94} $P_{n}^{(-\mu,\mu)}(z)$; in terms of
the latter, one has
\begin{eqnarray}
P_{n}^{\mu}(z)=\frac{n!}{\Gamma(n-\mu+1)}
\left(\frac{z+1}{z-1}\right)^{\mu/2}P_{n}^{(-\mu,\mu)}(z).
\label{2.3}
\end{eqnarray}
Combining this with the Rodrigues-type formula for the Jacobi
polynomials, which is
\begin{equation}
P_{n}^{(\alpha,\beta)}(z)=\frac{1}{2^{n}n!}
(z-1)^{-\alpha}(z+1)^{-\beta}
\frac{\mathrm{d}^{n}}{\mathrm{d}z^{n}}
\left[(z-1)^{n+\alpha}(z+1)^{n+\beta}\right],
\label{2.4}
\end{equation}
leads to the following Rodrigues-type representation of
$P_{n}^{\mu}(z)$:
\begin{equation}
P_{n}^{\mu}(z)=\frac{1}{2^{n}\Gamma(n-\mu+1)}
\left(\frac{z-1}{z+1}\right)^{\mu/2}
\frac{\mathrm{d}^{n}}{\mathrm{d}z^{n}}
\left[(z-1)^{n-\mu}(z+1)^{n+\mu}\right].
\label{2.5}
\end{equation}

The function $P_{n}^{\mu}(z)$ is single-valued in the complex
$z$-plane with a cut along the real axis from $-\infty$ to $1$. For a
negative integer degree, the associated Legendre function of the
first kind is defined through the relationship
\begin{equation}
P_{-n-1}^{\mu}(z)=P_{n}^{\mu}(z).
\label{2.6}
\end{equation}

For the sake of later applications, it is convenient to rewrite Eq.\
(\ref{2.1}) as
\begin{equation}
P_{n}^{\mu}(z)=(\pm)^{n}\frac{\Gamma(n\mp\mu+1)}{\Gamma(n-\mu+1)}
\left(\frac{z+1}{z-1}\right)^{\mu/2}
\sum_{k=0}^{n}(\pm)^{k}\frac{(k+n)!}{k!(n-k)!\Gamma(k\mp\mu+1)}
\left(\frac{z\mp1}{2}\right)^{k}
\label{2.7}
\end{equation}
and Eq.\ (\ref{2.2}) as
\begin{eqnarray}
P_{n}^{\mu}(z) &=& 
n!\Gamma(n+\mu+1)\left(\frac{z+1}{z-1}\right)^{\mu/2}
\left(\frac{z\pm1}{2}\right)^{n}
\nonumber \\
&& \times\sum_{k=0}^{n}
\frac{1}{k!(n-k)!\Gamma(k\mp\mu+1)\Gamma(n-k\pm\mu+1)}
\left(\frac{z\mp1}{z\pm1}\right)^{k}.
\label{2.8}
\end{eqnarray}
From either of Eqs.\ (\ref{2.1}), (\ref{2.2}), (\ref{2.7}) or
(\ref{2.8}), the function $P_{n}^{\mu}(z)$ is seen to possess the
reflection property
\begin{equation}
P_{n}^{\mu}(z)
=(-)^{n}\frac{\Gamma(n+\mu+1)}{\Gamma(n-\mu+1)}P_{n}^{-\mu}(-z).
\label{2.9}
\end{equation}

If $\mu=\pm m$, then it holds that
\begin{equation}
P_{n}^{m}(z)\equiv0
\qquad (m>n),
\label{2.10}
\end{equation}
\begin{equation}
P_{n}^{\pm m}(-z)=(-)^{n}P_{n}^{\pm m}(z)
\qquad (0\leqslant m\leqslant n),
\label{2.11}
\end{equation}
\begin{equation}
P_{n}^{-m}(z)=\frac{(n-m)!}{(n+m)!}P_{n}^{m}(z)
\qquad (0\leqslant m\leqslant n).
\label{2.12}
\end{equation}
If $0\leqslant m\leqslant n$, from Eqs.\ (\ref{2.7}) and (\ref{2.8})
we have
\begin{eqnarray}
P_{n}^{m}(z) &=& \left(\frac{z^{2}-1}{4}\right)^{m/2}
\sum_{k=0}^{n-m}\frac{(k+n+m)!}{k!(k+m)!(n-m-k)!}
\left(\frac{z-1}{2}\right)^{k}
\qquad (0\leqslant m\leqslant n),
\nonumber \\
\label{2.13}
\end{eqnarray}
\begin{eqnarray}
P_{n}^{m}(z) &=& (-)^{n}\frac{(n+m)!}{(n-m)!}
\left(\frac{z+1}{z-1}\right)^{m/2}
\sum_{k=0}^{n}(-)^{k}\frac{(k+n)!}{k!(k+m)!(n-k)!}
\left(\frac{z+1}{2}\right)^{k}
\nonumber \\
&& (0\leqslant m\leqslant n),
\label{2.14}
\end{eqnarray}
\begin{eqnarray}
P_{n}^{m}(z) &=& n!(n+m)!\left(\frac{z\mp1}{z\pm1}\right)^{m/2}
\left(\frac{z\pm1}{2}\right)^{n}
\nonumber \\
&& \times\sum_{k=0}^{n-m}\frac{1}{k!(k+m)!(n-k)!(n-m-k)!}
\left(\frac{z\mp1}{z\pm1}\right)^{k}
\qquad (0\leqslant m\leqslant n),
\label{2.15}
\end{eqnarray}
and also
\begin{equation}
P_{n}^{-m}(z)=\left(\frac{z-1}{z+1}\right)^{m/2}
\sum_{k=0}^{n}\frac{(k+n)!}{k!(k+m)!(n-k)!}
\left(\frac{z-1}{2}\right)^{k}
\qquad (0\leqslant m\leqslant n),
\label{2.16}
\end{equation}
\begin{eqnarray}
P_{n}^{-m}(z) &=& (-)^{n+m}\frac{(n-m)!}{(n+m)!}
\left(\frac{z^{2}-1}{4}\right)^{m/2}
\nonumber \\
&& \times\sum_{k=0}^{n-m}(-)^{k}\frac{(k+n+m)!}{k!(k+m)!(n-m-k)!}
\left(\frac{z+1}{2}\right)^{k}
\qquad (0\leqslant m\leqslant n),
\label{2.17}
\end{eqnarray}
\begin{eqnarray}
P_{n}^{-m}(z) &=& n!(n-m)!\left(\frac{z\mp1}{z\pm1}\right)^{m/2}
\left(\frac{z\pm1}{2}\right)^{n}
\nonumber \\
&& \times\sum_{k=0}^{n-m}\frac{1}{k!(k+m)!(n-k)!(n-m-k)!}
\left(\frac{z\mp1}{z\pm1}\right)^{k}
\qquad (0\leqslant m\leqslant n).
\label{2.18}
\end{eqnarray}
For $m>n$, the same equations yield
\begin{eqnarray}
P_{n}^{-m}(z) &=& \left(\frac{z-1}{z+1}\right)^{m/2}
\sum_{k=0}^{n}\frac{(k+n)!}{k!(k+m)!(n-k)!}
\left(\frac{z-1}{2}\right)^{k}
\qquad (m>n),
\label{2.19}
\end{eqnarray}
\begin{eqnarray}
P_{n}^{-m}(z) &=& \frac{1}{(n+m)!(m-n-1)!}
\left(\frac{z-1}{z+1}\right)^{m/2}
\sum_{k=0}^{n}\frac{(k+n)!(m-k-1)!}{k!(n-k)!}
\left(\frac{z+1}{2}\right)^{k}
\nonumber \\
&& (m>n),
\label{2.20}
\end{eqnarray}
\begin{eqnarray}
P_{n}^{-m}(z) &=& \frac{n!}{(m-n-1)!}
\left(\frac{z-1}{z+1}\right)^{m/2}
\left(\frac{z+1}{2}\right)^{n}
\sum_{k=0}^{n}(-)^{k}\frac{(k+m-n-1)!}{k!(k+m)!(n-k)!}
\left(\frac{z-1}{z+1}\right)^{k}
\nonumber \\
&& (m>n),
\label{2.21}
\end{eqnarray}
and
\begin{eqnarray}
P_{n}^{-m}(z) &=& (-)^{n}\frac{n!}{(m-n-1)!}
\left(\frac{z-1}{z+1}\right)^{m/2}
\left(\frac{z-1}{2}\right)^{n}
\nonumber \\
&& \times\sum_{k=0}^{n}(-)^{k}\frac{(m-k-1)!}{k!(n-k)!(n+m-k)!}
\left(\frac{z+1}{z-1}\right)^{k}
\qquad (m>n).
\label{2.22}
\end{eqnarray}
In deriving Eqs.\ (\ref{2.20}) to (\ref{2.22}), the formula
\begin{equation}
\lim_{\mu\to-m}\frac{\Gamma(n+\mu+1)}{\Gamma(n'+\mu+1)}
=(-)^{n+n'}\frac{(m-n'-1)!}{(m-n-1)!}
\qquad (m>n,n')
\label{2.23}
\end{equation}
appears to be useful.

Following Hobson \cite{Hobs31}, on that part of the cut for which
$z=x$ ($-1\leqslant x\leqslant1$) we define the associated
Legendre function of the first kind of non-negative degree as
\begin{equation}
P_{n}^{\mu}(x)=\frac{1}{2}
\left[\mathrm{e}^{\mathrm{i}\pi\mu/2}P_{n}^{\mu}(x+\mathrm{i}0)
+\mathrm{e}^{-\mathrm{i}\pi\mu/2}P_{n}^{\mu}(x-\mathrm{i}0)\right]
=\mathrm{e}^{\pm\mathrm{i}\pi\mu/2}P_{n}^{\mu}(x\pm\mathrm{i}0).
\label{2.24}
\end{equation}
\subsection{The associated Legendre function of the second kind of
integer degree}
\label{II.2}
The associated Legendre function of the second kind of a non-negative
integer degree $n$ and an arbitrary complex order $\mu$ may be
defined by either of the two expressions:
\begin{equation}
Q_{n}^{\mu}(z)=\frac{\pi}{2}
\frac{\mathrm{e}^{\mathrm{i}\pi\mu}}{\sin(\pi\mu)}
\left[P_{n}^{\mu}(z)
-\frac{\Gamma(n+\mu+1)}{\Gamma(n-\mu+1)}P_{n}^{-\mu}(z)\right],
\label{2.25}
\end{equation}
\begin{equation}
Q_{n}^{\mu}(z)=\frac{\pi}{2}
\frac{\mathrm{e}^{\mathrm{i}\pi\mu}}{\sin(\pi\mu)}
\left[P_{n}^{\mu}(z)-(-)^{n}P_{n}^{\mu}(-z)\right],
\label{2.26}
\end{equation}
equivalence of which follows immediately from Eq.\ (\ref{2.9}). An
extension of the definition to negative integer degrees is made via
the formula
\begin{equation}
Q_{-n-1}^{\mu}(z)=-Q_{n}^{\mu}(z)
+\pi\frac{\mathrm{e}^{\mathrm{i}\pi\mu}}
{\sin(\pi\mu)}P_{n}^{\mu}(z).
\label{2.27}
\end{equation}
One obtains
\begin{equation}
Q_{-n-1}^{\mu}(z)=\frac{\pi}{2}
\frac{\mathrm{e}^{\mathrm{i}\pi\mu}}{\sin(\pi\mu)}
\left[P_{n}^{\mu}(z)
+\frac{\Gamma(n+\mu+1)}{\Gamma(n-\mu+1)}P_{n}^{-\mu}(z)\right]
\label{2.28}
\end{equation}
and
\begin{equation}
Q_{-n-1}^{\mu}(z)=\frac{\pi}{2}
\frac{\mathrm{e}^{\mathrm{i}\pi\mu}}{\sin(\pi\mu)}
\left[P_{n}^{\mu}(z)+(-)^{n}P_{n}^{\mu}(-z)\right].
\label{2.29}
\end{equation}
From Eq.\ (\ref{2.25}) it is seen that
\begin{equation}
Q_{n}^{-\mu}(z)=\mathrm{e}^{-\mathrm{i}2\pi\mu}
\frac{\Gamma(n-\mu+1)}{\Gamma(n+\mu+1)}Q_{n}^{\mu}(z).
\label{2.30}
\end{equation}
It should be observed that the function $Q_{-n-1}^{\mu}(z)$ does
not exist if $\mu-n$ is a negative integer or zero. The functions
$Q_{n}^{\mu}(z)$ and $Q_{-n-1}^{\mu}(z)$ (provided the latter
exists, see above) are single-valued in the complex plane with a cut
along the real axis from $-\infty$ to $1$.

For $0\leqslant m\leqslant n$, it is
known\footnote[2]{\label{FOOT2}~In \cite[pages 81, 82, and
85]{Robi58} expressions for $Q_{n}^{m}(z)$, equivalent to these
following from our Eqs.\ (\ref{2.31}) to (\ref{2.33}), have been
given in somewhat different forms (in this connection, cf footnote
\ref{FOOT1} on page \pageref{FOOT1}). We have modified Robin's
formulas to make them concurrent with the notation used in the
present paper.} (see
\cite[pages 81, 82, and 85]{Robi58}) that
\begin{eqnarray}
Q_{n}^{m}(z)=\frac{1}{2}P_{n}^{m}(z)\ln\frac{z+1}{z-1}
-W_{n-1}^{m}(z)
\qquad (0\leqslant m\leqslant n),
\label{2.31}
\end{eqnarray}
with
\begin{eqnarray}
W_{n-1}^{m}(z) &=& 
\frac{1}{2}[\psi(n+m+1)+\psi(n-m+1)]P_{n}^{m}(z)
\nonumber \\
&& -\frac{(-)^{m}}{2}\left(\frac{z+1}{z-1}\right)^{m/2}
\sum_{k=0}^{m-1}(-)^{k}\frac{(k+n)!(m-k-1)!}{k!(n-k)!}
\left(\frac{z-1}{2}\right)^{k}
\nonumber \\
&& -\frac{1}{2}\frac{(n+m)!}{(n-m)!}
\left(\frac{z-1}{z+1}\right)^{m/2}
\sum_{k=0}^{n}\frac{(k+n)!\psi(k+m+1)}{k!(k+m)!(n-k)!}
\left(\frac{z-1}{2}\right)^{k}
\nonumber \\
&& -\frac{1}{2}\left(\frac{z^{2}-1}{4}\right)^{m/2}
\sum_{k=0}^{n-m}\frac{(k+n+m)!\psi(k+1)}{k!(k+m)!(n-m-k)!}
\left(\frac{z-1}{2}\right)^{k}
\qquad (0\leqslant m\leqslant n)
\nonumber \\
\label{2.32}
\end{eqnarray}
or
\begin{eqnarray}
W_{n-1}^{m}(z) &=&
-\frac{1}{2}[\psi(n+m+1)+\psi(n-m+1)]P_{n}^{m}(z)
\nonumber \\
&& +\frac{(-)^{n+m}}{2}\left(\frac{z-1}{z+1}\right)^{m/2}
\sum_{k=0}^{m-1}\frac{(k+n)!(m-k-1)!}{k!(n-k)!}
\left(\frac{z+1}{2}\right)^{k}
\nonumber \\
&& +\frac{(-)^{n}}{2}\frac{(n+m)!}{(n-m)!}
\left(\frac{z+1}{z-1}\right)^{m/2}
\sum_{k=0}^{n}(-)^{k}\frac{(k+n)!\psi(k+m+1)}{k!(k+m)!(n-k)!}
\left(\frac{z+1}{2}\right)^{k}
\nonumber \\
&& +\frac{(-)^{n+m}}{2}\left(\frac{z^{2}-1}{4}\right)^{m/2}
\sum_{k=0}^{n-m}(-)^{k}\frac{(k+n+m)!\psi(k+1)}{k!(k+m)!(n-m-k)!}
\left(\frac{z+1}{2}\right)^{k}
\nonumber \\
&& (0\leqslant m\leqslant n).
\label{2.33}
\end{eqnarray}
Another expression for $W_{n-1}^{m}(z)$ may be deduced from the
findings of Brown \cite{Brow95}; it is
\begin{eqnarray}
W_{n-1}^{m}(z)
&=& -\frac{(-)^{m}}{2}\left(\frac{z+1}{z-1}\right)^{m/2}
\sum_{k=0}^{m-1}(-)^{k}\frac{(k+n)!(m-k-1)!}{k!(n-k)!}
\left(\frac{z-1}{2}\right)^{k}
\nonumber \\
&& +\frac{(-)^{n+m}}{2}\left(\frac{z-1}{z+1}\right)^{m/2}
\sum_{k=0}^{m-1}\frac{(k+n)!(m-k-1)!}{k!(n-k)!}
\left(\frac{z+1}{2}\right)^{k}
\nonumber \\
&& -\frac{1}{2}\left(\frac{z^{2}-1}{4}\right)^{m/2}
\sum_{k=0}^{n-m}\frac{(k+n+m)!\psi(k+1)}{k!(k+m)!(n-m-k)!}
\left(\frac{z-1}{2}\right)^{k}
\nonumber \\
&& +\frac{(-)^{n+m}}{2}\left(\frac{z^{2}-1}{4}\right)^{m/2}
\sum_{k=0}^{n-m}(-)^{k}\frac{(k+n+m)!\psi(k+1)}{k!(k+m)!(n-m-k)!}
\left(\frac{z+1}{2}\right)^{k}
\nonumber \\
&& (0\leqslant m\leqslant n).
\label{2.34}
\end{eqnarray}
Furthermore, Snow \cite[pages 55 and 56]{Snow52} gave a
representation of $W_{n-1}^{m}(z)$ which may be easily shown to be
equivalent to
\begin{eqnarray}
W_{n-1}^{m}(z) &=& \pm\frac{1}{2}n!(n+m)!
\left(\frac{z\mp1}{z\pm1}\right)^{m/2}\left(\frac{z\pm1}{2}\right)^{n}
\sum_{k=0}^{n-m}\frac{1}{k!(k+m)!(n-k)!(n-m-k)!}
\nonumber \\
&& \quad\times[\psi(n-m-k+1)+\psi(n-k+1)-\psi(k+m+1)-\psi(k+1)]
\left(\frac{z\mp1}{z\pm1}\right)^{k}
\nonumber \\
&& \pm\frac{1}{2}n!(n+m)!\left(\frac{z\pm1}{z\mp1}\right)^{m/2}
\left(\frac{z\mp1}{2}\right)^{n}
\nonumber \\
&& \quad\times\sum_{k=1}^{m}(-)^{k}\frac{(k-1)!}{(k+n)!(k+n-m)!(m-k)!}
\left(\frac{z\mp1}{z\pm1}\right)^{k}
\nonumber \\
&& \mp\frac{(-)^{m}}{2}n!(n+m)!\left(\frac{z\pm1}{z\mp1}\right)^{m/2}
\left(\frac{z\pm1}{2}\right)^{n}
\nonumber \\
&& \quad\times\sum_{k=0}^{m-1}(-)^{k}
\frac{(m-k-1)!}{k!(n-k)!(n+m-k)!}
\left(\frac{z\mp1}{z\pm1}\right)^{k}
\qquad (0\leqslant m\leqslant n).
\label{2.35}
\end{eqnarray}

Exploiting Eq.\ (\ref{2.31}) and either of Eqs.\ (\ref{2.32}) to
(\ref{2.34}), $Q_{n}^{-m}(z)$ may be evaluated using
\begin{equation}
Q_{n}^{-m}(z)=\frac{(n-m)!}{(n+m)!}Q_{n}^{m}(z)
\qquad (0\leqslant m\leqslant n),
\label{2.36}
\end{equation}
which is the direct consequence of Eq.\ (\ref{2.30}). The functions
$Q_{-n-1}^{\pm m}(z)$ with $0\leqslant m\leqslant n$ do not
exist.

For $m>n$, it holds that \cite[Eq.\ (63) on page 35]{Robi58}
\begin{equation}
Q_{n}^{m}(z)=\frac{(-)^{m}}{2}(n+m)!(m-n-1)!
\left[P_{n}^{-m}(-z)-(-)^{n}P_{n}^{-m}(z)\right]
\qquad (m>n).
\label{2.37}
\end{equation}
The function $Q_{n}^{-m}(z)$ with $m>n$ does not exist.

On the cut $x\in[-1,1]$, again following Hobson \cite{Hobs31}, we
define the associated Legendre function of the second kind of
non-negative and negative integer degrees as
\begin{equation}
Q_{n}^{\mu}(x)=\frac{1}{2}\mathrm{e}^{-\mathrm{i}\pi\mu}
\left[\mathrm{e}^{-\mathrm{i}\pi\mu/2}Q_{n}^{\mu}(x+\mathrm{i}0)
+\mathrm{e}^{\mathrm{i}\pi\mu/2}Q_{n}^{\mu}(x-\mathrm{i}0)\right]
\label{2.38}
\end{equation}
and
\begin{equation}
Q_{-n-1}^{\mu}(x)=\frac{1}{2}\mathrm{e}^{-\mathrm{i}\pi\mu}
\left[\mathrm{e}^{-\mathrm{i}\pi\mu/2}Q_{-n-1}^{\mu}(x+\mathrm{i}0)
+\mathrm{e}^{\mathrm{i}\pi\mu/2}
Q_{-n-1}^{\mu}(x-\mathrm{i}0)\right],
\label{2.39}
\end{equation}
respectively.
%
%
\section{The derivatives $\partial P_{n}^{\mu}(z)/\partial\mu$
and $\partial P_{n}^{\mu}(x)/\partial\mu$}
\label{III}
\setcounter{equation}{0}
\subsection{General considerations}
\label{III.1}
We begin with the observation that after differentiating Eqs.\
(\ref{2.6}) and (\ref{2.9}) with respect to $\mu$, one arrives at the
following two general relations involving the derivative in question:
\begin{equation}
\frac{\partial P_{-n-1}^{\mu}(z)}{\partial\mu}
=\frac{\partial P_{n}^{\mu}(z)}{\partial\mu},
\label{3.1}
\end{equation}
\begin{eqnarray}
\frac{\partial P_{n}^{-\mu}(z)}{\partial\mu}
&=& -[\psi(n+\mu+1)+\psi(n-\mu+1)]P_{n}^{-\mu}(z)
+(-)^{n}\frac{\Gamma(n-\mu+1)}{\Gamma(n+\mu+1)}
\frac{\partial P_{n}^{\mu}(-z)}{\partial\mu}.
\nonumber \\
\label{3.2}
\end{eqnarray}
Because of Eq.\ (\ref{3.1}), in the remainder of the paper we shall
focus on evaluation of $\partial P_{n}^{\mu}(z)/\partial\mu$.

We proceed to the construction of several explicit representations of
$\partial P_{n}^{\mu}(z)/\partial\mu$.

From Eq.\ (\ref{2.5}), we have the Rodrigues-type formula
\begin{eqnarray}
\frac{\partial P_{n}^{\mu}(z)}{\partial\mu}
&=& -\frac{1}{2}P_{n}^{\mu}(z)\ln\frac{z+1}{z-1}
+\psi(n-\mu+1)P_{n}^{\mu}(z)
\nonumber \\
&& +\frac{1}{2^{n}\Gamma(n-\mu+1)}
\left(\frac{z-1}{z+1}\right)^{\mu/2}
\frac{\mathrm{d}^{n}}{\mathrm{d}z^{n}}
\left[(z-1)^{n-\mu}(z+1)^{n+\mu}\ln\frac{z+1}{z-1}\right].
\label{3.3}
\end{eqnarray}
This may be rewritten as
\begin{equation}
\frac{\partial P_{n}^{\mu}(z)}{\partial\mu}
=\frac{1}{2}P_{n}^{\mu}(z)\ln\frac{z+1}{z-1}+U_{n}^{\mu}(z)
\label{3.4}
\end{equation}
with
\begin{eqnarray}
U_{n}^{\mu}(z) &=& -P_{n}^{\mu}(z)\ln\frac{z+1}{z-1}
+\psi(n-\mu+1)P_{n}^{\mu}(z)
\nonumber \\
&& +\frac{1}{2^{n}\Gamma(n-\mu+1)}
\left(\frac{z-1}{z+1}\right)^{\mu/2}
\frac{\mathrm{d}^{n}}{\mathrm{d}z^{n}}
\left[(z-1)^{n-\mu}(z+1)^{n+\mu}\ln\frac{z+1}{z-1}\right].
\label{3.5}
\end{eqnarray}

Choosing in Eq.\ (\ref{2.7}) the upper set of signs and
differentiating with respect to $\mu$ yields $\partial
P_{n}^{\mu}(z)/\partial\mu$ in the form (\ref{3.4}) with
\begin{equation}
U_{n}^{\mu}(z)=\left(\frac{z+1}{z-1}\right)^{\mu/2}
\sum_{k=0}^{n}\frac{(k+n)!\psi(k-\mu+1)}{k!(n-k)!\Gamma(k-\mu+1)}
\left(\frac{z-1}{2}\right)^{k}.
\label{3.6}
\end{equation}
An alternative representation of $U_{n}^{\mu}(z)$ is obtained if in
Eq.\ (\ref{2.7}) one chooses the lower set of signs and then
differentiates with respect to $\mu$; this results in
\begin{eqnarray}
U_{n}^{\mu}(z) &=& [\psi(n+\mu+1)+\psi(n-\mu+1)]P_{n}^{\mu}(z)
-(-)^{n}\frac{\Gamma(n+\mu+1)}{\Gamma(n-\mu+1)}
\left(\frac{z+1}{z-1}\right)^{\mu/2}
\nonumber \\
&& \times\sum_{k=0}^{n}(-)^{k}
\frac{(k+n)!\psi(k+\mu+1)}{k!(n-k)!\Gamma(k+\mu+1)}
\left(\frac{z+1}{2}\right)^{k}.
\label{3.7}
\end{eqnarray}
Playing with Eq.\ (\ref{3.7}) with the aid of Eq.\ (\ref{2.7}) (with
the lower signs chosen) and of the following easy to prove
identities:
\begin{equation}
\frac{\Gamma(n+\mu+1)}{\Gamma(n'+\mu+1)}
=(-)^{n+n'}\frac{\Gamma(-n'-\mu)}{\Gamma(-n-\mu)}
\label{3.8}
\end{equation}
and
\begin{equation}
\psi(n+\mu+1)-\psi(n'+\mu+1)=\psi(-n-\mu)-\psi(-n'-\mu)
\label{3.9}
\end{equation}
leads to
\begin{eqnarray}
U_{n}^{\mu}(z) &=& [\psi(n-\mu+1)+\psi(-n-\mu)]P_{n}^{\mu}(z)
-\frac{1}{\Gamma(n-\mu+1)\Gamma(-n-\mu)}
\left(\frac{z+1}{z-1}\right)^{\mu/2}
\nonumber \\
&& \times\sum_{k=0}^{n}
\frac{(k+n)!\Gamma(-k-\mu)\psi(-k-\mu)}{k!(n-k)!}
\left(\frac{z+1}{2}\right)^{k}.
\label{3.10}
\end{eqnarray}
Furthermore, differentiating Eq.\ (\ref{2.8}) with respect to $\mu$
gives
\begin{eqnarray}
U_{n}^{\mu}(z) &=& \psi(n+\mu+1)P_{n}^{\mu}(z)
\pm n!\Gamma(n+\mu+1)\left(\frac{z+1}{z-1}\right)^{\mu/2}
\left(\frac{z\pm1}{2}\right)^{n}
\nonumber \\
&& \times\sum_{k=0}^{n}\frac{\psi(k\mp\mu+1)-\psi(n-k\pm\mu+1)}
{k!(n-k)!\Gamma(k\mp\mu+1)\Gamma(n-k\pm\mu+1)}
\left(\frac{z\mp1}{z\pm1}\right)^{k}.
\label{3.11}
\end{eqnarray}
Manipulations with Eq.\ (\ref{3.11}), similar to those which have led
us from Eq.\ (\ref{3.7}) to Eq.\ (\ref{3.10}), give
\begin{eqnarray}
U_{n}^{\mu}(z) &=& \psi(-n-\mu)P_{n}^{\mu}(z)
+\frac{n!}{\Gamma(-n-\mu)}\left(\frac{z+1}{z-1}\right)^{\mu/2}
\left(\frac{z+1}{2}\right)^{n}
\nonumber \\
&& \times\sum_{k=0}^{n}(-)^{k}\frac{\Gamma(k-n-\mu)}
{k!(n-k)!\Gamma(k-\mu+1)}[\psi(k-\mu+1)-\psi(k-n-\mu)]
\left(\frac{z-1}{z+1}\right)^{k}
\nonumber \\
\label{3.12}
\end{eqnarray}
or
\begin{eqnarray}
U_{n}^{\mu}(z) &=& \psi(-n-\mu)P_{n}^{\mu}(z)
+(-)^{n}\frac{n!}{\Gamma(-n-\mu)}
\left(\frac{z+1}{z-1}\right)^{\mu/2}\left(\frac{z-1}{2}\right)^{n}
\nonumber \\
&& \times\sum_{k=0}^{n}(-)^{k}\frac{\Gamma(-k-\mu)}
{k!(n-k)!\Gamma(n-k-\mu+1)}[\psi(n-k-\mu+1)-\psi(-k-\mu)]
\left(\frac{z+1}{z-1}\right)^{k},
\nonumber \\
\label{3.13}
\end{eqnarray}
according to whether the upper or the lower signs are chosen in Eq.\
(\ref{3.11}).

Next, let us look at Eqs.\ (\ref{2.7}) and (\ref{3.6}) which give
the functions $P_{n}^{\mu}(z)$ and $U_{n}^{\mu}(z)$,
respectively. It is seen that both these functions are of the form
\mbox{$[(z+1)/(z-1)]^{\mu/2}$} times a polynomial in $z-1$ of degree
$n$ [except for the case of $\mu=m>n$, when $P_{n}^{m}(z)$ vanishes
identically, see Eq.\ (\ref{2.10})]. Consequently, if $\mu\neq m$, it
must be possible to represent $U_{n}^{\mu}(z)$ in the form of a
linear combination of the functions $P_{k}^{\mu}(z)$ with degrees not
exceeding $n$:
\begin{equation}
U_{n}^{\mu}(z)=\sum_{k=0}^{n}c_{nk}^{\mu}P_{k}^{\mu}(z)
\label{3.14}
\end{equation}
[in the case of $\mu=m$ one may seek $U_{n}^{\mu}(z)$ in the form of
an expansion analogous to (\ref{3.14}) but with $P_{k}^{\mu}(z)$
replaced by $P_{k}^{-\mu}(-z)$; cf also Eq.\ (\ref{3.34})]. Comparing
coefficients at $(z-1)^{n}[(z+1)/(z-1)]^{\mu/2}$ in Eqs.\ (\ref{2.7})
and (\ref{3.6}), we see that the coefficient $c_{nn}^{\mu}$ in the
above expansion is
\begin{equation}
c_{nn}^{\mu}=\psi(n-\mu+1).
\label{3.15}
\end{equation}
The reasoning leading to determination of the coefficients
$c_{nk}^{\mu}$ with $0\leqslant k\leqslant n-1$ is as follows.

If we differentiate the well-known Legendre identity
\cite{Magn66,Grad94}
\begin{equation}
\left[\frac{\mathrm{d}}{\mathrm{d}z}(1-z^{2})
\frac{\mathrm{d}}{\mathrm{d}z}+n(n+1)-\frac{\mu^{2}}{1-z^{2}}\right]
P_{n}^{\mu}(z)=0
\label{3.16}
\end{equation}
with respect to $\mu$, this yields
\begin{equation}
\left[\frac{\mathrm{d}}{\mathrm{d}z}(1-z^{2})
\frac{\mathrm{d}}{\mathrm{d}z}+n(n+1)-\frac{\mu^{2}}{1-z^{2}}\right]
\frac{\partial P_{n}^{\mu}(z)}{\partial\mu}
=\frac{2\mu}{1-z^{2}}P_{n}^{\mu}(z).
\label{3.17}
\end{equation}
On exploiting the identity (\ref{3.16}), we also obtain
\begin{equation}
\left[\frac{\mathrm{d}}{\mathrm{d}z}(1-z^{2})
\frac{\mathrm{d}}{\mathrm{d}z}+n(n+1)-\frac{\mu^{2}}{1-z^{2}}\right]
\frac{1}{2}P_{n}^{\mu}(z)\ln\frac{z+1}{z-1}
=2\frac{\mathrm{d}P_{n}^{\mu}(z)}{\mathrm{d}z}.
\label{3.18}
\end{equation}
By subtracting Eq.\ (\ref{3.18}) from Eq.\ (\ref{3.17}), we find that
the function $U_{n}^{\mu}(z)$ satisfies the inhomogeneous
differential equation
\begin{equation}
\left[\frac{\mathrm{d}}{\mathrm{d}z}(1-z^{2})
\frac{\mathrm{d}}{\mathrm{d}z}+n(n+1)-\frac{\mu^{2}}{1-z^{2}}\right]
U_{n}^{\mu}(z)=-2\frac{\mathrm{d}P_{n}^{\mu}(z)}{\mathrm{d}z}
+\frac{2\mu}{1-z^{2}}P_{n}^{\mu}(z).
\label{3.19}
\end{equation}

For a while, let us focus on the expression on the right-hand side of
Eq.\ (\ref{3.19}). We shall show how it may be developed into a
finite sum of the functions $P_{k}^{\mu}(z)$ with $0\leqslant
k\leqslant n-1$. To begin with, we observe that it holds that
\begin{eqnarray}
-2\frac{\mathrm{d}P_{n}^{\mu}(z)}{\mathrm{d}z}
+\frac{2\mu}{1-z^{2}}P_{n}^{\mu}(z)
&=&
\left[(z-1)\frac{\mathrm{d}P_{n}^{\mu}(z)}{\mathrm{d}z}
+\frac{\mu}{z+1}P_{n}^{\mu}(z)-nP_{n}^{\mu}(z)\right]
\nonumber \\
&& -\left[(z+1)\frac{\mathrm{d}P_{n}^{\mu}(z)}{\mathrm{d}z}
+\frac{\mu}{z-1}P_{n}^{\mu}(z)-nP_{n}^{\mu}(z)\right].
\label{3.20}
\end{eqnarray}
If use is made of the well-known relation \cite{Magn66,Grad94}
\begin{equation}
(z^{2}-1)\frac{\mathrm{d}P_{n}^{\mu}(z)}{\mathrm{d}z}
=nzP_{n}^{\mu}(z)-(n+\mu)P_{n-1}^{\mu}(z),
\label{3.21}
\end{equation}
the expressions in the square brackets appearing on the right-hand
side of Eq.\ (\ref{3.20}) may be written as
\begin{equation}
(z\mp1)\frac{\mathrm{d}P_{n}^{\mu}(z)}{\mathrm{d}z}
+\frac{\mu}{z\pm1}P_{n}^{\mu}(z)-nP_{n}^{\mu}(z)
=\mp\frac{(n\mp\mu)P_{n}^{\mu}(z)
\pm(n+\mu)P_{n-1}^{\mu}(z)}{z\pm1}.
\label{3.22}
\end{equation}
Consider next the recurrence relation \cite{Magn66,Grad94}
\begin{equation}
(2k+1)zP_{k}^{\mu}(z)=(k-\mu+1)P_{k+1}^{\mu}(z)
+(k+\mu)P_{k-1}^{\mu}(z).
\label{3.23}
\end{equation}
Evidently, it may rewritten as
\begin{eqnarray}
(2k+1)(z+1)P_{k}^{\mu}(z)
&=& \left[(k-\mu+1)P_{k+1}^{\mu}(z)+(k+\mu+1)P_{k}^{\mu}(z)\right]
\nonumber \\
&& +\left[(k-\mu)P_{k}^{\mu}(z)+(k+\mu)P_{k-1}^{\mu}(z)\right].
\label{3.24}
\end{eqnarray}
Multiplying both sides of Eq.\ (\ref{3.24}) by $(-)^{k+n}$, summing
over $k$ from $k=0$ to $k=n-1$, utilizing then the property [cf Eq.\
(\ref{2.6})]
\begin{equation}
P_{0}^{\mu}(z)=P_{-1}^{\mu}(z)
\label{3.25}
\end{equation}
and dividing the result by $z+1$, we find
\begin{equation}
\sum_{k=0}^{n-1}(-)^{k+n}(2k+1)P_{k}^{\mu}(z)
=-\frac{(n-\mu)P_{n}^{\mu}(z)+(n+\mu)P_{n-1}^{\mu}(z)}{z+1}.
\label{3.26}
\end{equation}
If on both sides of Eq.\ (\ref{3.26}) use is made of Eq.\
(\ref{2.9}), the former becomes
\begin{equation}
\sum_{k=0}^{n-1}(2k+1)\frac{\Gamma(k+\mu+1)}{\Gamma(k-\mu+1)}
P_{k}^{-\mu}(-z)=-\frac{\Gamma(n+\mu+1)}{\Gamma(n-\mu)}
\frac{P_{n}^{-\mu}(-z)-P_{n-1}^{-\mu}(-z)}{z+1}.
\label{3.27}
\end{equation}
From this, after replacing $\mu$ with $-\mu$ and $z$ with $-z$, we
infer
\begin{equation}
\sum_{k=0}^{n-1}(2k+1)\frac{\Gamma(k-\mu+1)}{\Gamma(k+\mu+1)}
P_{k}^{\mu}(z)=\frac{\Gamma(n-\mu+1)}{\Gamma(n+\mu)}
\frac{P_{n}^{\mu}(z)-P_{n-1}^{\mu}(z)}{z-1}.
\label{3.28}
\end{equation}
Hence, it follows that
\begin{eqnarray}
\frac{\Gamma(n+\mu+1)}{\Gamma(n-\mu+1)}
\sum_{k=0}^{n-1}(2k+1)
\frac{\Gamma(k-\mu+1)}{\Gamma(k+\mu+1)}P_{k}^{\mu}(z)
=\frac{(n+\mu)P_{n}^{\mu}(z)-(n+\mu)P_{n-1}^{\mu}(z)}{z-1}.
\label{3.29}
\end{eqnarray}
Combining Eqs.\ (\ref{3.20}), (\ref{3.22}), (\ref{3.28}), and
(\ref{3.29}), we arrive at the sought expansion
\begin{eqnarray}
-2\frac{\mathrm{d}P_{n}^{\mu}(z)}{\mathrm{d}z}
&+& \frac{2\mu}{1-z^{2}}P_{n}^{\mu}(z)
\nonumber \\
&=& \sum_{k=0}^{n-1}(-)^{k+n}(2k+1)
\left[1-(-)^{k+n}\frac{\Gamma(n+\mu+1)\Gamma(k-\mu+1)}
{\Gamma(n-\mu+1)\Gamma(k+\mu+1)}\right]P_{k}^{\mu}(z),
\nonumber \\
\label{3.30}
\end{eqnarray}
which, apart from being useful in the present context, seems to be
also interesting for its own sake.

We return to the problem of determination of the coefficients
$c_{nk}^{\mu}$. In virtue of the result in Eq.\ (\ref{3.30}), we
rewrite the differential relation (\ref{3.19}) as
\begin{eqnarray}
&& \left[\frac{\mathrm{d}}{\mathrm{d}z}(1-z^{2})
\frac{\mathrm{d}}{\mathrm{d}z} 
+n(n+1)-\frac{\mu^{2}}{1-z^{2}}\right]U_{n}^{\mu}(z)
\nonumber \\
&& \qquad\qquad =\:\sum_{k=0}^{n-1}(-)^{k+n}(2k+1)
\left[1-(-)^{k+n}\frac{\Gamma(n+\mu+1)\Gamma(k-\mu+1)}
{\Gamma(n-\mu+1)\Gamma(k+\mu+1)}\right]P_{k}^{\mu}(z).
\label{3.31}
\end{eqnarray}
On substituting the expansion (\ref{3.14}) into the left-hand side of
Eq.\ (\ref{3.31}), after equating coefficients at $P_{k}^{\mu}(z)$
on both sides of the resulting identity, we obtain
\begin{eqnarray}
c_{nk}^{\mu} &=& (-)^{k+n}\frac{2k+1}{(n-k)(k+n+1)}
\left[1-(-)^{k+n}\frac{\Gamma(n+\mu+1)\Gamma(k-\mu+1)}
{\Gamma(n-\mu+1)\Gamma(k+\mu+1)}\right]
\nonumber \\
&& (0\leqslant k\leqslant n-1).
\label{3.32}
\end{eqnarray}
Plugging Eqs.\ (\ref{3.15}) and (\ref{3.32}) into Eq.\ (\ref{3.14}),
we arrive at the following representation of $U_{n}^{\mu}(z)$:
\begin{eqnarray}
U_{n}^{\mu}(z) &=& \psi(n-\mu+1)P_{n}^{\mu}(z)
+\sum_{k=0}^{n-1}(-)^{k+n}\frac{2k+1}{(n-k)(k+n+1)}
\nonumber \\
&& \times\left[1-(-)^{k+n}\frac{\Gamma(n+\mu+1)\Gamma(k-\mu+1)}
{\Gamma(n-\mu+1)\Gamma(k+\mu+1)}\right]P_{k}^{\mu}(z).
\label{3.33}
\end{eqnarray}
With the use of the property (\ref{2.9}), the above result may be
transformed into
\begin{eqnarray}
U_{n}^{\mu}(z)
&=& \psi(n-\mu+1)P_{n}^{\mu}(z)
\nonumber \\
&& +\sum_{k=0}^{n-1}(-)^{k+n}\frac{2k+1}{(n-k)(k+n+1)}
\left[P_{k}^{\mu}(z)-(-)^{n}\frac{\Gamma(n+\mu+1)}{\Gamma(n-\mu+1)}
P_{k}^{-\mu}(-z)\right].
\nonumber \\
\label{3.34}
\end{eqnarray}
The last expression is of value when $0\leqslant\mu=m\leqslant n$.
%
%
\subsection{Evaluation of $[\partial
P_{n}^{\mu}(z)/\partial\mu]_{\mu=m}$}
\label{III.2}
\subsubsection{The case of $0\leqslant m\leqslant n$}
\label{III.2.1}
For $\mu=m$, the fundamental equation (\ref{3.4}) becomes
\begin{equation}
\frac{\partial P_{n}^{\mu}(z)}{\partial\mu}\bigg|_{\mu=m}
=\frac{1}{2}P_{n}^{m}(z)\ln\frac{z+1}{z-1}+U_{n}^{m}(z)
\qquad (0\leqslant m\leqslant n).
\label{3.35}
\end{equation}
It follows from Eq.\ (\ref{3.5}) that the Rodrigues-type
representation of $U_{n}^{m}(z)$ is
\begin{eqnarray}
U_{n}^{m}(z) &=& -P_{n}^{m}(z)\ln\frac{z+1}{z-1}
+\psi(n-m+1)P_{n}^{m}(z)
\nonumber \\
&& +\frac{1}{2^{n}(n-m)!}\left(\frac{z-1}{z+1}\right)^{m/2}
\frac{\mathrm{d}^{n}}{\mathrm{d}z^{n}}
\left[(z-1)^{n-m}(z+1)^{n+m}\ln\frac{z+1}{z-1}\right]
\nonumber \\
&& (0\leqslant m\leqslant n).
\label{3.36}
\end{eqnarray}

If in Eq.\ (\ref{3.6}) use is made of the well-known limiting
property
\begin{equation}
\lim_{\mu\to m}\frac{\psi(k-\mu+1)}{\Gamma(k-\mu+1)}
=(-)^{k+m}(m-k-1)!
\qquad (m>k)
\label{3.37}
\end{equation}
and if, whenever necessary, the gamma function is replaced by the
factorial, after some straightforward manipulations we obtain the
following representation of $U_{n}^{m}(z)$:
\begin{eqnarray}
U_{n}^{m}(z) &=& (-)^{m}\left(\frac{z+1}{z-1}\right)^{m/2}
\sum_{k=0}^{m-1}(-)^{k}\frac{(k+n)!(m-k-1)!}{k!(n-k)!}
\left(\frac{z-1}{2}\right)^{k}
\nonumber \\
&& +\left(\frac{z^{2}-1}{4}\right)^{m/2}\sum_{k=0}^{n-m}
\frac{(k+n+m)!\psi(k+1)}{k!(k+m)!(n-m-k)!}
\left(\frac{z-1}{2}\right)^{k}
\qquad (0\leqslant m\leqslant n).
\nonumber \\
\label{3.38}
\end{eqnarray}
Similarly, from Eq.\ (\ref{3.7}) it follows that
\begin{eqnarray}
U_{n}^{m}(z) &=& [\psi(n+m+1)+\psi(n-m+1)]P_{n}^{m}(z)
\nonumber \\
&& -(-)^{n}\frac{(n+m)!}{(n-m)!}\left(\frac{z+1}{z-1}\right)^{m/2}
\sum_{k=0}^{n}(-)^{k}
\frac{(k+n)!\psi(k+m+1)}{k!(k+m)!(n-k)!}
\left(\frac{z+1}{2}\right)^{k}
\nonumber \\
&& (0\leqslant m\leqslant n).
\label{3.39}
\end{eqnarray}
Furthermore, from Eq.\ (\ref{3.11}) one finds
\begin{eqnarray}
U_{n}^{m}(z) &=& \psi(n+m+1)P_{n}^{m}(z)
+(-)^{m}n!(n+m)!\left(\frac{z+1}{z-1}\right)^{m/2}
\left(\frac{z+1}{2}\right)^{n}
\nonumber \\
&& \quad\times\sum_{k=0}^{m-1}(-)^{k}
\frac{(m-k-1)!}{k!(n-k)!(n+m-k)!}\left(\frac{z-1}{z+1}\right)^{k}
\nonumber \\
&& -n!(n+m)!\left(\frac{z-1}{z+1}\right)^{m/2}
\left(\frac{z+1}{2}\right)^{n}
\nonumber \\
&& \quad\times\sum_{k=0}^{n-m}
\frac{\psi(n-k+1)-\psi(k+1)}{k!(k+m)!(n-k)!(n-m-k)!}
\left(\frac{z-1}{z+1}\right)^{k}
\qquad (0\leqslant m\leqslant n)
\label{3.40}
\end{eqnarray}
and
\begin{eqnarray}
U_{n}^{m}(z) &=& \psi(n+m+1)P_{n}^{m}(z)
+n!(n+m)!\left(\frac{z-1}{z+1}\right)^{m/2}
\left(\frac{z+1}{2}\right)^{n}
\nonumber \\
&& \quad\times\sum_{k=1}^{m}(-)^{k}
\frac{(k-1)!}{(k+n)!(k+n-m)!(m-k)!}\left(\frac{z+1}{z-1}\right)^{k}
\nonumber \\
&& +n!(n+m)!\left(\frac{z+1}{z-1}\right)^{m/2}
\left(\frac{z-1}{2}\right)^{n}
\nonumber \\
&& \quad\times\sum_{k=0}^{n-m}
\frac{\psi(n-m-k+1)-\psi(k+m+1)}{k!(k+m)!(n-k)!(n-m-k)!}
\left(\frac{z+1}{z-1}\right)^{k}
\qquad (0\leqslant m\leqslant n).
\nonumber \\
&&
\label{3.41}
\end{eqnarray}

Finally, from either of Eqs.\ (\ref{3.33}) or (\ref{3.34}), with the
aid of Eqs.\ (\ref{2.13}) to (\ref{2.15}), one deduces that
\begin{eqnarray}
U_{n}^{m}(z) &=& \psi(n-m+1)P_{n}^{m}(z)-\frac{(n+m)!}{(n-m)!}
\sum_{k=0}^{n-1}(-)^{k}\frac{2k+1}{(n-k)(k+n+1)}P_{k}^{-m}(-z)
\nonumber \\
&& +\sum_{k=0}^{n-m-1}(-)^{k+n+m}\frac{2k+2m+1}{(n-m-k)(k+n+m+1)}
P_{k+m}^{m}(z)
\qquad (0\leqslant m\leqslant n)
\nonumber \\
&&
\label{3.42}
\end{eqnarray}
or equivalently
\begin{eqnarray}
U_{n}^{m}(z) &=& \psi(n-m+1)P_{n}^{m}(z)-\frac{(n+m)!}{(n-m)!}
\sum_{k=0}^{m-1}(-)^{k}\frac{2k+1}{(n-k)(k+n+1)}P_{k}^{-m}(-z)
\nonumber \\
&& +\sum_{k=0}^{n-m-1}(-)^{k+n+m}\frac{2k+2m+1}{(n-m-k)(k+n+m+1)}
\nonumber \\
&& \quad 
\times\left[1-(-)^{k+n+m}\frac{k!(n+m)!}{(k+2m)!(n-m)!}\right]
P_{k+m}^{m}(z)
\qquad (0\leqslant m\leqslant n).
\label{3.43}
\end{eqnarray}
\subsubsection{A closer look at the case of $m=0$}
\label{III.2.2}
Setting $\mu=0$ in Eq.\ (\ref{3.17}) results in
\begin{equation}
\left[\frac{\mathrm{d}}{\mathrm{d}z}(1-z^{2})
\frac{\mathrm{d}}{\mathrm{d}z}+n(n+1)\right]
\frac{\partial P_{n}^{\mu}(z)}{\partial\mu}\bigg|_{\mu=0}=0,
\label{3.44}
\end{equation}
i.e., $[\partial P_{n}^{\mu}(z)/\partial\mu]_{\mu=0}$ solves the
same differential equation as the Legendre polynomial $P_{n}(z)\equiv
P_{n}^{0}(z)$ does. Consequently, it must be possible to write it
in the form of a linear combination of $P_{n}(z)$ and the Legendre
function of the second kind $Q_{n}(z)\equiv Q_{n}^{0}(z)$.

To find that combination, at first we observe that for $m=0$ Eq.\
(\ref{3.35}) becomes
\begin{equation}
\frac{\partial P_{n}^{\mu}(z)}{\partial\mu}\bigg|_{\mu=0}
=\frac{1}{2}P_{n}(z)\ln\frac{z+1}{z-1}+U_{n}^{0}(z).
\label{3.45}
\end{equation}
From Eq.\ (\ref{3.43}) we see that $U_{n}^{0}(z)$ is a polynomial in
$z$ of degree $n$, which may be written as
\begin{equation}
U_{n}^{0}(z)=\psi(n+1)P_{n}(z)-W_{n-1}(z).
\label{3.46}
\end{equation}
Here $W_{n-1}(z)$ is a polynomial in $z$ (of degree $n-1$) given by
\begin{equation}
W_{n-1}(z)=2\sum_{k=0}^{n-1}\frac{1-(-)^{k+n}}{2}
\frac{2k+1}{(n-k)(k+n+1)}P_{k}(z)
\label{3.47}
\end{equation}
or equivalently
\begin{equation}
W_{n-1}(z)=\sum_{k=0}^{\mathrm{int}[(n-1)/2]}
\frac{2n-4k-1}{(n-k)(2k+1)}P_{n-2k-1}(z).
\label{3.48}
\end{equation}
A glance at Eq.\ (\ref{3.48}) reveals that $W_{n-1}(z)$ is the
well-known Christoffel's polynomial (cf, e.g., \cite[page
153]{Erde53}), in terms of which one has \cite{Magn66,Grad94}
\begin{equation}
Q_{n}(z)=\frac{1}{2}P_{n}(z)\ln\frac{z+1}{z-1}-W_{n-1}(z).
\label{3.49}
\end{equation}
On combining Eqs.\ (\ref{3.45}), (\ref{3.46}), and (\ref{3.49}), we
obtain the sought relationship\footnote[3]{\label{FOOT3}~The
counterpart relationship (see \cite[page 178]{Magn66})
\begin{displaymath}
\frac{\partial Q_{n}^{\mu}(z)}{\partial\mu}\bigg|_{\mu=0}
=[\mathrm{i}\pi+\psi(n+1)]Q_{n}(z)
\end{displaymath}
follows straightforwardly from Eq.\ (\ref{2.30}) after
differentiating the latter with respect to $\mu$ and setting then
$\mu=0$.} [cf Eq.\ (\ref{1.1})]
\begin{equation}
\frac{\partial P_{n}^{\mu}(z)}{\partial\mu}\bigg|_{\mu=0}
=\psi(n+1)P_{n}(z)+Q_{n}(z).
\label{3.50}
\end{equation}
The same result may be obtained if one couples Eq.\ (\ref{3.2}) with
Eq.\ (\ref{4.1}), both particularized to the case of $m=0$.
%
%
\subsubsection{The case of $m>n$}
\label{III.2.3}
In this case, the function $P_{n}^{m}(z)$ vanishes identically [cf
Eq.\ (\ref{2.10})], so that Eq.\ (\ref{3.4}) becomes
\begin{equation}
\frac{\partial P_{n}^{\mu}(z)}{\partial\mu}\bigg|_{\mu=m}
=U_{n}^{m}(z)
\qquad (m>n).
\label{3.51}
\end{equation}
With the help of Eq.\ (\ref{3.37}), from Eq.\ (\ref{3.6}) we find
\begin{eqnarray}
U_{n}^{m}(z) &=& (-)^{m}\left(\frac{z+1}{z-1}\right)^{m/2}
\sum_{k=0}^{n}(-)^{k}\frac{(k+n)!(m-k-1)!}{k!(n-k)!}
\left(\frac{z-1}{2}\right)^{k}
\qquad (m>n).
\nonumber \\
\label{3.52}
\end{eqnarray}
Comparison of Eq.\ (\ref{3.52}) with Eq.\ (\ref{2.20}) reveals the
relationship\footnote[4]{\label{FOOT4}~The result (\ref{3.53}) may be
also deduced from Eq.\ (\ref{3.34}). Indeed, it is seen that in the
case considered here the summand in the latter equation vanishes
identically and only the first term on the right-hand side survives.
Further, in virtue of Eqs.\ (\ref{2.9}) and (\ref{3.37}) (with $k$
replaced by $n$), we have
\begin{displaymath}
\lim_{\mu\to m}\psi(n-\mu+1)P_{n}^{\mu}(z)
=(-)^{m}(n+m)!(m-n-1)!P_{n}^{-m}(-z)
\qquad (m>n),
\end{displaymath}
hence, Eq.\ (\ref{3.53}) follows. Still another way of arriving at
Eq.\ (\ref{3.53}) is to use Eq.\ (\ref{3.5}).}
\begin{equation}
U_{n}^{m}(z)=(-)^{m}(n+m)!(m-n-1)!P_{n}^{-m}(-z)
\qquad (m>n).
\label{3.53}
\end{equation}
Consequently, from Eqs.\ (\ref{3.51}) and (\ref{3.53}) one obtains
\begin{equation}
\frac{\partial P_{n}^{\mu}(z)}{\partial\mu}\bigg|_{\mu=m}
=(-)^{m}(n+m)!(m-n-1)!P_{n}^{-m}(-z)
\qquad (m>n).
\label{3.54}
\end{equation}
%
%
\subsection{Evaluation of $[\partial
P_{n}^{\mu}(z)/\partial\mu]_{\mu=-m}$}
\label{III.3}
\subsubsection{The case of $0\leqslant m\leqslant n$}
\label{III.3.1}
For $\mu=-m$, Eq.\ (\ref{3.4}) becomes
\begin{equation}
\frac{\partial P_{n}^{\mu}(z)}{\partial\mu}\bigg|_{\mu=-m}
=\frac{1}{2}P_{n}^{-m}(z)\ln\frac{z+1}{z-1}
+U_{n}^{-m}(z)
\qquad (0\leqslant m\leqslant n).
\label{3.55}
\end{equation}
Since it holds that
\begin{equation}
\frac{\partial P_{n}^{\mu}(z)}{\partial\mu}\bigg|_{\mu=-m}
=-\frac{\partial P_{n}^{-\mu}(z)}{\partial\mu}\bigg|_{\mu=m},
\label{3.56}
\end{equation}
from Eqs.\ (\ref{3.55}), (\ref{3.2}), and (\ref{3.35}) it may be
inferred that
\begin{eqnarray}
U_{n}^{-m}(z) &=& (-)^{n+1}\frac{(n-m)!}{(n+m)!}U_{n}^{m}(-z)
+[\psi(n+m+1)+\psi(n-m+1)]P_{n}^{-m}(z)
\nonumber \\
&& (0\leqslant m\leqslant n).
\label{3.57}
\end{eqnarray}
One may couple this formula with Eqs.\ (\ref{3.38}) to (\ref{3.43})
and exploit, whenever necessary, Eqs.\ (\ref{2.11}) and (\ref{2.12})
to obtain various explicit representations of $U_{n}^{-m}(z)$ with
$0\leqslant m\leqslant n$.
\subsubsection{The case of $m>n$}
\label{III.3.2}
From Eq.\ (\ref{3.4}) we find
\begin{equation}
\frac{\partial P_{n}^{\mu}(z)}{\partial\mu}\bigg|_{\mu=-m}
=\frac{1}{2}P_{n}^{-m}(z)\ln\frac{z+1}{z-1}+U_{n}^{-m}(z)
\qquad (m>n).
\label{3.58}
\end{equation}
Equations (\ref{3.6}), (\ref{3.10}), (\ref{3.12}), and (\ref{3.13})
supply us with the following alternative representations of
$U_{n}^{-m}(z)$:
\begin{equation}
U_{n}^{-m}(z)=\left(\frac{z-1}{z+1}\right)^{m/2}
\sum_{k=0}^{n}\frac{(k+n)!\psi(k+m+1)}{k!(k+m)!(n-k)!}
\left(\frac{z-1}{2}\right)^{k}
\qquad (m>n),
\label{3.59}
\end{equation}
\begin{eqnarray}
U_{n}^{-m}(z) &=& [\psi(n+m+1)+\psi(m-n)]P_{n}^{-m}(z)
-\frac{1}{(n+m)!(m-n-1)!}\left(\frac{z-1}{z+1}\right)^{m/2}
\nonumber \\
&& \times\sum_{k=0}^{n}\frac{(k+n)!(m-k-1)!\psi(m-k)}{k!(n-k)!}
\left(\frac{z+1}{2}\right)^{k}
\qquad (m>n),
\label{3.60}
\end{eqnarray}
\begin{eqnarray}
U_{n}^{-m}(z) &=& \psi(m-n)P_{n}^{-m}(z)
+\frac{n!}{(m-n-1)!}\left(\frac{z-1}{z+1}\right)^{m/2}
\left(\frac{z+1}{2}\right)^{n}
\nonumber \\
&& \times\sum_{k=0}^{n}(-)^{k}\frac{(k+m-n-1)!}{k!(k+m)!(n-k)!}
[\psi(k+m+1)-\psi(k+m-n)]\left(\frac{z-1}{z+1}\right)^{k}
\nonumber \\
&& (m>n),
\label{3.61}
\end{eqnarray}
\begin{eqnarray}
U_{n}^{-m}(z) &=& \psi(m-n)P_{n}^{-m}(z)
+(-)^{n}\frac{n!}{(m-n-1)!}\left(\frac{z-1}{z+1}\right)^{m/2}
\left(\frac{z-1}{2}\right)^{n}
\nonumber \\
&& \times\sum_{k=0}^{n}(-)^{k}\frac{(m-k-1)!}{k!(n-k)!(n+m-k)!}
[\psi(n+m-k+1)-\psi(m-k)]\left(\frac{z+1}{z-1}\right)^{k}
\nonumber \\
&& (m>n).
\label{3.62}
\end{eqnarray}
Moreover, with the aid of Eq.\ (\ref{2.23}), from Eq.\ (\ref{3.33})
we obtain
\begin{eqnarray}
U_{n}^{-m}(z) &=& \psi(n+m+1)P_{n}^{-m}(z)
\nonumber \\
&& +\sum_{k=0}^{n-1}(-)^{k+n}\frac{2k+1}{(n-k)(k+n+1)}
\left[1-\frac{(k+m)!(m-k-1)!}{(n+m)!(m-n-1)!}\right]
P_{k}^{-m}(z)
\nonumber \\
&& (m>n).
\label{3.63}
\end{eqnarray}
%
%
\subsection{The derivative $\partial P_{n}^{\mu}(x)/\partial\mu$}
\label{III.4}
Since [cf Eq.\ (\ref{2.24})] the definition of $P_{n}^{\mu}(x)$
with $-1\leqslant x\leqslant1$ differs from that of
$P_{n}^{\mu}(z)$ with $z\in\mathbb{C}\setminus[-1,1]$, the
derivative $\partial P_{n}^{\mu}(x)/\partial\mu$ requires to be
considered separately.

Differentiating Eq.\ (\ref{2.24}) with respect to $\mu$, making use
of Eq.\ (\ref{3.35}) and exploiting the fact that
\begin{equation}
x+1\pm\mathrm{i}0=1+x
\qquad x-1\pm\mathrm{i}0=\mathrm{e}^{\pm\mathrm{i}\pi}(1-x)
\qquad(-1\leqslant x\leqslant1)
\label{3.64}
\end{equation}
results in
\begin{equation}
\frac{\partial P_{n}^{\mu}(x)}{\partial\mu}
=\frac{1}{2}P_{n}^{\mu}(x)\ln\frac{1+x}{1-x}+U_{n}^{\mu}(x),
\label{3.65}
\end{equation}
where, in analogy with Eq.\ (\ref{2.24}), we define
\begin{eqnarray}
U_{n}^{\mu}(x)=\frac{1}{2}
\left[\mathrm{e}^{\mathrm{i}\pi\mu/2}U_{n}^{\mu}(x+\mathrm{i}0)
+\mathrm{e}^{-\mathrm{i}\pi\mu/2}U_{n}^{\mu}(x-\mathrm{i}0)\right]
=\mathrm{e}^{\pm\mathrm{i}\pi\mu/2}U_{n}^{\mu}(x\pm\mathrm{i}0).
\label{3.66}
\end{eqnarray}
Various representations of $U_{n}^{\mu}(x)$ and $U_{n}^{\pm m}(x)$
arise if one couples Eq.\ (\ref{3.66}) with the formulas for
$U_{n}^{\mu}(z)$ and $U_{n}^{\pm m}(z)$ derived in sections
\ref{III.1} to \ref{III.3} (now particularized to the case of
$z=x\pm\mathrm{i}0$).
%
%
\section{Applications}
\label{IV}
\setcounter{equation}{0}
\subsection{Construction of the associated Legendre function of the
second kind of integer degree and order}
\label{IV.1}
As the first application of the results obtained in the preceding
section, below we shall construct several explicit representations of
the associated Legendre function of the second kind of integer degree
and order.
\subsubsection{The functions $Q_{n}^{\pm m}(z)$ in the case of
$0\leqslant m\leqslant n$}
\label{IV.1.1}
If in Eq.\ (\ref{2.26}) we let $\mu$ approach $\pm m$, after applying
the l'Hospital rule to the right-hand side, we obtain
\begin{equation}
Q_{n}^{\pm m}(z)
=\frac{1}{2}\frac{\partial P_{n}^{\mu}(z)}
{\partial\mu}\bigg|_{\mu=\pm m}
-\frac{(-)^{n}}{2}\frac{\partial P_{n}^{\mu}(-z)}
{\partial\mu}\bigg|_{\mu=\pm m}
\qquad (0\leqslant m\leqslant n).
\label{4.1}
\end{equation}
Combining this with Eq.\ (\ref{3.3}) gives the following
Rodrigues-type representation of $Q_{n}^{\pm m}(z)$ with
$0\leqslant m\leqslant n$:
\begin{eqnarray}
Q_{n}^{\pm m}(z) 
&=& -\frac{1}{2}P_{n}^{\pm m}(z)\ln\frac{z+1}{z-1}
\nonumber \\
&& +\frac{1}{2^{n+1}(n\mp m)!}\left(\frac{z-1}{z+1}\right)^{m/2}
\frac{\mathrm{d}^{n}}{\mathrm{d}z^{n}}
\left[(z-1)^{n-m}(z+1)^{n+m}\ln\frac{z+1}{z-1}\right]
\nonumber \\
&& +\frac{1}{2^{n+1}(n\mp m)!}\left(\frac{z+1}{z-1}\right)^{m/2}
\frac{\mathrm{d}^{n}}{\mathrm{d}z^{n}}
\left[(z-1)^{n+m}(z+1)^{n-m}\ln\frac{z+1}{z-1}\right]
\nonumber \\
&& (0\leqslant m\leqslant n).
\label{4.2}
\end{eqnarray}
For $m=0$, Eq.\ (\ref{4.2}) reduces to the well-known (e.g.,
\cite[Eq.\ (8.836.1)]{Grad94}) formula
\begin{equation}
Q_{n}(z)=-\frac{1}{2}P_{n}(z)\ln\frac{z+1}{z-1}
+\frac{1}{2^{n}n!}\frac{\mathrm{d}^{n}}{\mathrm{d}z^{n}}
\left[(z^{2}-1)^{n}\ln\frac{z+1}{z-1}\right].
\label{4.3}
\end{equation}

Upon making use of Eqs.\ (\ref{3.4}) and (\ref{2.11}), Eq.\
(\ref{4.1}) implies that
\begin{equation}
Q_{n}^{\pm m}(z)=\frac{1}{2}P_{n}^{\pm m}(z)\ln\frac{z+1}{z-1}
-W_{n-1}^{\pm m}(z)
\qquad (0\leqslant m\leqslant n),
\label{4.4}
\end{equation}
with
\begin{equation}
W_{n-1}^{\pm m}(z)=\frac{1}{2}\left[(-)^{n}U_{n}^{\pm m}(-z)
-U_{n}^{\pm m}(z)\right]
\qquad (0\leqslant m\leqslant n).
\label{4.5}
\end{equation}
From Eq.\ (\ref{4.5}) it is seen that
\begin{equation}
W_{n-1}^{\pm m}(-z)=(-)^{n+1}W_{n-1}^{\pm m}(z)
\qquad (0\leqslant m\leqslant n),
\label{4.6}
\end{equation}
i.e., the parities of $W_{n-1}^{\pm m}(z)$ are opposite to those of
$P_{n}^{\pm m}(z)$. In view of the symmetry relation
\begin{equation}
W_{n-1}^{-m}(z)=\frac{(n-m)!}{(n+m)!}W_{n-1}^{m}(z)
\qquad (0\leqslant m\leqslant n),
\label{4.7}
\end{equation}
which is the consequence of Eqs.\ (\ref{4.5}), (\ref{3.57}), and
(\ref{2.11}), henceforth we shall discuss $W_{n-1}^{m}(z)$ only.

A number of explicit representations of $W_{n-1}^{m}(z)$ may be
obtained if one couples Eq.\ (\ref{4.5}) (with the plus signs chosen)
with Eqs.\ (\ref{3.38}) to (\ref{3.43}). For instance, if
$U_{n}^{m}(-z)$ is calculated using Eq.\ (\ref{3.39}), while
$U_{n}^{m}(z)$ is taken in the form (\ref{3.38}), then Eq.\
(\ref{2.32}) is recovered. Conversely, if $U_{n}^{m}(-z)$ is obtained
from Eq.\ (\ref{3.38}) and $U_{n}^{m}(z)$ from Eq.\ (\ref{3.39}),
this results in Eq.\ (\ref{2.33}). If in Eq.\ (\ref{4.5}) both
$U_{n}^{m}(-z)$ and $U_{n}^{m}(z)$ are obtained from Eq.\
(\ref{3.38}), one recovers the Brown's formula (\ref{2.34}), while if
both $U_{n}^{m}(-z)$ and $U_{n}^{m}(z)$ are derived from Eq.\
(\ref{3.39}), one arrives at
\begin{eqnarray}
W_{n-1}^{m}(z) &=& -\frac{1}{2}\frac{(n+m)!}{(n-m)!}
\left(\frac{z-1}{z+1}\right)^{m/2}
\sum_{k=0}^{n}\frac{(k+n)!\psi(k+m+1)}{k!(k+m)!(n-k)!}
\left(\frac{z-1}{2}\right)^{k}
\nonumber \\
&& +\frac{(-)^{n}}{2}\frac{(n+m)!}{(n-m)!}
\left(\frac{z+1}{z-1}\right)^{m/2}
\sum_{k=0}^{n}(-)^{k}\frac{(k+n)!\psi(k+m+1)}{k!(k+m)!(n-k)!}
\left(\frac{z+1}{2}\right)^{k}
\nonumber \\
&& (0\leqslant m\leqslant n).
\label{4.8}
\end{eqnarray}
That of the two Snow's representations which corresponds to the
choice of the upper signs in Eq.\ (\ref{2.35}) follows if in Eq.\
(\ref{4.5}) $U_{n}^{m}(-z)$ is found from Eq.\ (\ref{3.41}) and
$U_{n}^{m}(z)$ from Eq.\ (\ref{3.40}); the other one is obtained if
the roles played by Eqs.\ (\ref{3.40}) and (\ref{3.41}) are
interchanged. 

An apparently new representation of $W_{n-1}^{m}(z)$ follows if in
Eq.\ (\ref{4.5}) both $U_{n}^{m}(-z)$ and $U_{n}^{m}(z)$ are found
from Eq.\ (\ref{3.33}). Proceeding in this way, with the aid of Eq.\
(\ref{2.11}), one finds
\begin{eqnarray}
W_{n-1}^{m}(z) &=& \frac{1}{2}\frac{(n+m)!}{(n-m)!}
\sum_{k=0}^{m-1}(-)^{k}\frac{2k+1}{(n-k)(k+n+1)}
\left[P_{k}^{-m}(-z)-(-)^{n}P_{k}^{-m}(z)\right]
\nonumber \\
&& +\sum_{k=0}^{n-m-1}\frac{1-(-)^{k+n+m}}{2}
\frac{2k+2m+1}{(n-m-k)(k+n+m+1)}
\nonumber \\
&& \quad\times\left[1-(-)^{k+n+m}\frac{k!(n+m)!}{(k+2m)!(n-m)!}\right]
P_{k+m}^{m}(z)
\qquad (0\leqslant m\leqslant n).
\label{4.9}
\end{eqnarray}
It is seen that the summand in the second sum on the right-hand side
of the above equation vanishes if $k+n+m$ is even and therefore Eq.\
(\ref{4.9}) may be rewritten more suitably as
\begin{eqnarray}
W_{n-1}^{m}(z) &=& \frac{1}{2}\frac{(n+m)!}{(n-m)!}
\sum_{k=0}^{m-1}(-)^{k}\frac{2k+1}{(n-k)(k+n+1)}
\left[P_{k}^{-m}(-z)-(-)^{n}P_{k}^{-m}(z)\right]
\nonumber \\
&& +\sum_{k=0}^{n-m-1}\frac{1-(-)^{k+n+m}}{2}
\frac{2k+2m+1}{(n-m-k)(k+n+m+1)}
\nonumber \\
&& \quad\times\left[1+\frac{k!(n+m)!}{(k+2m)!(n-m)!}\right]
P_{k+m}^{m}(z)
\qquad (0\leqslant m\leqslant n).
\label{4.10}
\end{eqnarray}
Manipulating with the second sum, Eq.\ (\ref{4.10}) may be cast into
\begin{eqnarray}
W_{n-1}^{m}(z) &=& \frac{1}{2}\frac{(n+m)!}{(n-m)!}
\sum_{k=0}^{m-1}(-)^{k}\frac{2k+1}{(n-k)(k+n+1)}
\left[P_{k}^{-m}(-z)-(-)^{n}P_{k}^{-m}(z)\right]
\nonumber \\
&& +\frac{1}{2}\sum_{k=0}^{\mathrm{int}[(n-m-1)/2]}
\frac{2n-4k-1}{(n-k)(2k+1)}
\nonumber \\
&& \quad
\times\left[1+\frac{(n+m)!(n-m-2k-1)!}{(n-m)!(n+m-2k-1)!}\right]
P_{n-2k-1}^{m}(z)
\qquad (0\leqslant m\leqslant n).
\label{4.11}
\end{eqnarray}
It is seen that in the case of $m=0$ Eqs.\ (\ref{4.10}) and
(\ref{4.11}) reduce to Eqs.\ (\ref{3.47}) and (\ref{3.48}),
respectively.
%
%
\subsubsection{The functions $Q_{n}^{m}(z)$ and $Q_{-n-1}^{m}(z)$
in the case of $m>n$}
\label{IV.1.2}
Proceeding as in section \ref{IV.1.1}, from Eq.\ (\ref{2.26}) we
obtain
\begin{equation}
Q_{n}^{m}(z)=\frac{1}{2}\frac{\partial P_{n}^{\mu}(z)}
{\partial\mu}\bigg|_{\mu=m}
-\frac{(-)^{n}}{2}\frac{\partial P_{n}^{\mu}(-z)}
{\partial\mu}\bigg|_{\mu=m}
\qquad (m>n).
\label{4.12}
\end{equation}
If use is made here of Eq.\ (\ref{3.54}), this yields Eq.\
(\ref{2.37}). Next, for the function $Q_{-n-1}^{m}(z)$, from Eq.\
(\ref{2.27}), with the aid of the l'Hospital rule, we have
\begin{equation}
Q_{-n-1}^{m}(z)=-Q_{n}^{m}(z)
+\frac{\partial P_{n}^{\mu}(z)}{\partial\mu}\bigg|_{\mu=m}
\qquad (m>n).
\label{4.13}
\end{equation}
From this, after merging with Eqs.\ (\ref{2.37}) and (\ref{3.54}),
we find
\begin{equation}
Q_{-n-1}^{m}(z)=\frac{(-)^{m}}{2}(n+m)!(m-n-1)!
\left[P_{n}^{-m}(-z)+(-)^{n}P_{n}^{-m}(z)\right]
\qquad (m>n).
\label{4.14}
\end{equation}
%
%
\subsubsection{The functions $Q_{n}^{\pm m}(x)$ and
$Q_{-n-1}^{m}(x)$} 
\label{IV.1.3}
Inserting Eq.\ (\ref{4.4}) into Eq.\ (\ref{2.38}), the latter being
particularized to $\mu=\pm m$, with the help of Eq.\ (\ref{2.24}) we
obtain
\begin{equation}
Q_{n}^{\pm m}(x)=\frac{1}{2}P_{n}^{\pm m}(x)\ln\frac{1+x}{1-x}
-W_{n-1}^{\pm m}(x)
\qquad (0\leqslant m\leqslant n),
\label{4.15}
\end{equation}
where we define
\begin{eqnarray}
W_{n-1}^{\pm m}(x) &=& \frac{(-)^{m}}{2}
\left[\mathrm{e}^{\mp\mathrm{i}\pi m/2}
W_{n-1}^{\pm m}(x+\mathrm{i}0)
+\mathrm{e}^{\pm \mathrm{i}\pi m/2}
W_{n-1}^{\pm m}(x-\mathrm{i}0)\right]
\qquad (0\leqslant m\leqslant n).
\nonumber \\
\label{4.16}
\end{eqnarray}
It follows from Eqs.\ (\ref{4.16}) and (\ref{4.6}) that the functions
$W_{n-1}^{\pm m}(x)$ possess the reflection property
\begin{equation}
W_{n-1}^{\pm m}(-x)=(-)^{n+m+1}W_{n-1}^{\pm m}(x)
\qquad (0\leqslant m\leqslant n),
\label{4.17}
\end{equation}
while if Eq.\ (\ref{4.16}) is coupled with Eq.\ (\ref{4.7}), this
results in
\begin{equation}
W_{n-1}^{-m}(x)=(-)^{m}\frac{(n-m)!}{(n+m)!}W_{n-1}^{m}(x)
\qquad (0\leqslant m\leqslant n).
\label{4.18}
\end{equation}
The latter relationship allows one to focus on $W_{n-1}^{m}(x)$
only.

Exploiting the expressions for $W_{n-1}^{m}(z)$ found in section
\ref{IV.1.1}, with no difficulty one may construct counterpart
formulas for $W_{n-1}^{m}(x)$. For instance, if Eq.\ (\ref{4.16}) is
combined with Eq.\ (\ref{4.11}) and use is made of Eq.\ (\ref{2.11}),
this yields
\begin{eqnarray}
W_{n-1}^{m}(x) &=& \frac{1}{2}\frac{(n+m)!}{(n-m)!}
\sum_{k=0}^{m-1}(-)^{k}\frac{2k+1}{(n-k)(k+n+1)}
\left[P_{k}^{-m}(-x)-(-)^{n+m}P_{k}^{-m}(x)\right]
\nonumber \\
&& +\frac{1}{2}\sum_{k=0}^{\mathrm{int}[(n-m-1)/2]}
\frac{2n-4k-1}{(n-k)(2k+1)}
\nonumber \\
&& \quad
\times\left[1+\frac{(n+m)!(n-m-2k-1)!}{(n-m)!(n+m-2k-1)!}\right]
P_{n-2k-1}^{m}(x)
\qquad (0\leqslant m\leqslant n).
\label{4.19}
\end{eqnarray}

If Eq.\ (\ref{2.37}) is used in Eq.\ (\ref{2.38}), one finds
\begin{equation}
Q_{n}^{m}(x)=\frac{(-)^{m}}{2}(n+m)!(m-n-1)!
\left[P_{n}^{-m}(-x)-(-)^{n+m}P_{n}^{-m}(x)\right]
\qquad (m>n),
\label{4.20}
\end{equation}
while if Eq.\ (\ref{4.14}) is coupled with Eq.\ (\ref{2.39}), this
gives
\begin{equation}
Q_{-n-1}^{m}(x)=\frac{(-)^{m}}{2}(n+m)!(m-n-1)!
\left[P_{n}^{-m}(-x)+(-)^{n+m}P_{n}^{-m}(x)\right]
\qquad (m>n).
\label{4.21}
\end{equation}
%
%
\subsection{Evaluation of
$[\partial^{2}P_{n}^{\mu}(z)/\partial\mu^{2}]_{\mu=m}$ for $m>n$}
\label{IV.2}
In this section, we shall show that if $m>n$, then the knowledge of
$[\partial P_{n}^{\mu}(-z)/\partial\mu]_{\mu=-m}$ allows one to
find $[\partial^{2}P_{n}^{\mu}(z)/\partial\mu^{2}]_{\mu=m}$.

We begin with the observation that from the identity
\begin{equation}
\psi(\zeta)=\psi(1-\zeta)-\cos(\pi\zeta)\Gamma(\zeta)\Gamma(1-\zeta),
\label{4.22}
\end{equation}
which may be easily deduced from elementary properties of the digamma
and gamma functions \cite{Erde53,Jahn60,Magn66,Grad94,Davi65}, one
finds that
\begin{equation}
\psi(n-\mu+1)=\psi(\mu-n)
+(-)^{n}\cos(\pi\mu)\Gamma(\mu-n)\Gamma(n-\mu+1).
\label{4.23}
\end{equation}
The identity (\ref{4.23}) allows one to rewrite the relation
(\ref{3.2}) in the form
\begin{eqnarray}
\frac{\partial P_{n}^{-\mu}(z)}{\partial\mu}
&=& -[\psi(n+\mu+1)+\psi(\mu-n)]P_{n}^{-\mu}(z)
\nonumber \\
&& +(-)^{n}\frac{\Gamma(n-\mu+1)}{\Gamma(n+\mu+1)}
\left[\frac{\partial P_{n}^{\mu}(-z)}{\partial\mu}
-\cos(\pi\mu)\Gamma(\mu-n)\Gamma(n+\mu+1)P_{n}^{-\mu}(z)\right],
\nonumber \\
\label{4.24}
\end{eqnarray}
hence, it follows that
\begin{eqnarray}
\frac{\partial P_{n}^{-\mu}(z)}{\partial\mu}\bigg|_{\mu=m}
&=& -[\psi(n+m+1)+\psi(m-n)]P_{n}^{-m}(z)
\nonumber \\
&& +\frac{(-)^{n}}{(n+m)!}\lim_{\mu\to m}\Bigg\{\Gamma(n-\mu+1)
\nonumber \\
&& \quad\times\Bigg[\frac{\partial P_{n}^{\mu}(-z)}{\partial\mu}
-\cos(\pi\mu)\Gamma(\mu-n)\Gamma(n+\mu+1)
P_{n}^{-\mu}(z)\Bigg]\Bigg\}
\nonumber \\
&& (m>n).
\label{4.25}
\end{eqnarray}
The limit in the above equation may be evaluated with the aid of the
l'Hospital rule [for this purpose, the property (\ref{3.37}) proves
helpful] and one finds
\begin{eqnarray}
\frac{\partial P_{n}^{-\mu}(z)}{\partial\mu}\bigg|_{\mu=m}
&=& -2[\psi(n+m+1)+\psi(m-n)]P_{n}^{-m}(z)
-\frac{\partial P_{n}^{-\mu}(z)}{\partial\mu}\bigg|_{\mu=m}
\nonumber \\
&& +\frac{(-)^{m}}{(n+m)!(m-n-1)!}
\frac{\partial^{2}P_{n}^{\mu}(-z)}{\partial\mu^{2}}\bigg|_{\mu=m}
\qquad (m>n).
\label{4.26}
\end{eqnarray}
If Eq.\ (\ref{4.26}) is solved for
$[\partial^{2}P_{n}^{\mu}(-z)/\partial\mu^{2}]_{\mu=m}$, then use
is made of Eq.\ (\ref{3.56}) and $z$ is replaced by $-z$, one
eventually arrives at
\begin{eqnarray}
\frac{\partial^{2}P_{n}^{\mu}(z)}{\partial\mu^{2}}\bigg|_{\mu=m}
&=& (-)^{m}2(n+m)!(m-n-1)!
\nonumber \\
&& \times\left\{[\psi(n+m+1)+\psi(m-n)]P_{n}^{-m}(-z)
-\frac{\partial P_{n}^{\mu}(-z)}{\partial\mu}\bigg|_{\mu=-m}\right\}
\nonumber \\
&& (m>n).
\label{4.27}
\end{eqnarray}
Several explicit representations of
$[\partial^{2}P_{n}^{\mu}(z)/\partial\mu^{2}]_{\mu=m}$ with $m>n$,
not listed here, follow if one combines Eq.\ (\ref{4.27}) with the
results of section \ref{III.3.2}.
%
%
\subsection{Evaluation of
$[\partial Q_{n}^{\mu}(z)/\partial\mu]_{\mu=m}$ and
$[\partial Q_{-n-1}^{\mu}(z)/\partial\mu]_{\mu=m}$ for $m>n$}
\label{IV.3}
The results of section \ref{IV.2} may be exploited to express the
derivatives $[\partial Q_{n}^{\mu}(z)/\partial\mu]_{\mu=m}$ and 
$[\partial Q_{-n-1}^{\mu}(z)/\partial\mu]_{\mu=m}$ for $m>n$ in
terms of $[\partial P_{n}^{\mu}(\pm z)/\partial\mu]_{\mu=-m}$. 

To show this for $[\partial Q_{n}^{\mu}(z)/\partial\mu]_{\mu=m}$,
we differentiate Eq.\ (\ref{2.26}) with respect to $\mu$, obtaining
\begin{eqnarray}
\frac{\partial Q_{n}^{\mu}(z)}{\partial\mu}
&=& \mathrm{i}\pi Q_{n}^{\mu}(z)+\frac{\pi}{\sin(\pi\mu)}
\Bigg[-\cos(\pi\mu)Q_{n}^{\mu}(z)
+\frac{1}{2}\mathrm{e}^{\mathrm{i}\pi\mu}
\frac{\partial P_{n}^{\mu}(z)}{\partial\mu}
-\frac{(-)^{n}}{2}\mathrm{e}^{\mathrm{i}\pi\mu}
\frac{\partial P_{n}^{\mu}(-z)}{\partial\mu}\Bigg].
\nonumber \\
\label{4.28}
\end{eqnarray}
In the limit $\mu\to m>n$, with the help of the l'Hospital rule and
of Eq.\ (\ref{4.12}), from Eq.\ (\ref{4.28}) we obtain
\begin{eqnarray}
\frac{\partial Q_{n}^{\mu}(z)}{\partial\mu}\bigg|_{\mu=m}
&=& 2\mathrm{i}\pi Q_{n}^{m}(z)
-\frac{\partial Q_{n}^{\mu}(z)}{\partial\mu}\bigg|_{\mu=m}
+\frac{1}{2}\frac{\partial^{2}P_{n}^{\mu}(z)}
{\partial\mu^{2}}\bigg|_{\mu=m}
-\frac{(-)^{n}}{2}\frac{\partial^{2}P_{n}^{\mu}(-z)}
{\partial\mu^{2}}\bigg|_{\mu=m}
\nonumber \\
&& (m>n),
\label{4.29}
\end{eqnarray}
hence, it follows that
\begin{eqnarray}
\frac{\partial Q_{n}^{\mu}(z)}{\partial\mu}\bigg|_{\mu=m}
&=& \mathrm{i}\pi Q_{n}^{m}(z)
+\frac{1}{4}\frac{\partial^{2}P_{n}^{\mu}(z)}
{\partial\mu^{2}}\bigg|_{\mu=m}
-\frac{(-)^{n}}{4}\frac{\partial^{2}P_{n}^{\mu}(-z)}
{\partial\mu^{2}}\bigg|_{\mu=m}
\qquad (m>n).
\nonumber \\
\label{4.30}
\end{eqnarray}
If the last two terms on the right-hand side of Eq.\ (\ref{4.30}) are
transformed with the aid of Eq.\ (\ref{4.27}) and use is made of Eq.\
(\ref{2.37}), this yields the sought relationship
\begin{eqnarray}
\frac{\partial Q_{n}^{\mu}(z)}{\partial\mu}\bigg|_{\mu=m}
&=& [\mathrm{i}\pi+\psi(n+m+1)+\psi(m-n)]Q_{n}^{m}(z)
\nonumber \\
&& -\frac{(-)^{m}}{2}(n+m)!(m-n-1)!
\left[\frac{\partial P_{n}^{\mu}(-z)}
{\partial\mu}\bigg|_{\mu=-m}
-(-)^{n}\frac{\partial P_{n}^{\mu}(z)}
{\partial\mu}\bigg|_{\mu=-m}\right]
\nonumber \\
&& (m>n).
\label{4.31}
\end{eqnarray}

If one starts with Eq.\ (\ref{2.29}) rather than with Eq.\
(\ref{2.26}), after movements almost identical to these presented
above one arrives at
\begin{eqnarray}
\frac{\partial Q_{-n-1}^{\mu}(z)}{\partial\mu}\bigg|_{\mu=m}
&=& [\mathrm{i}\pi+\psi(n+m+1)+\psi(m-n)]Q_{-n-1}^{m}(z)
\nonumber \\
&& -\frac{(-)^{m}}{2}(n+m)!(m-n-1)!
\left[\frac{\partial P_{n}^{\mu}(-z)}
{\partial\mu}\bigg|_{\mu=-m}
+(-)^{n}\frac{\partial P_{n}^{\mu}(z)}
{\partial\mu}\bigg|_{\mu=-m}\right]
\nonumber \\
&& (m>n).
\label{4.32}
\end{eqnarray}
%
%

%

\begin{thebibliography}{99}
\bibitem{Todh75}
   I.\ Todhunter,
   \emph{An Elementary Treatise on Laplace's Functions, Lam{\'e}'s
   Functions and Bessel's Functions\/}
   (Macmillan, London, 1875).
\bibitem{Ferr77}
   N.M.\ Ferrers,
   \emph{An Elementary Treatise on Spherical Harmonics\/}
   (Macmillan, London, 1877).
\bibitem{Neum78}
   F.\ Neumann,
   \emph{Beitr{\"a}ge zur Theorie der Kugelfunctionen\/}
   (Teubner, Leipzig, 1878).
\bibitem{Hein78}
   E.\ Heine,
   \emph{Handbuch der Kugelfunctionen\/}, vol.\ 1, 2nd ed.\
   (Reimer, Berlin, 1878).
\bibitem{Hein81}
   E.\ Heine,
   \emph{Handbuch der Kugelfunctionen\/}, vol.\ 2, 2nd ed.\
   (Reimer, Berlin, 1881).
\bibitem{Olbr87}
   R.\ Olbricht,
   Nova Acta Leop.\ Carol.\ Akad.\ 52 (1887) 1.
\bibitem{Byer93}
   W.E.\ Byerly,
   \emph{An Elementary Treatise on Fourier's Series and Spherical,
   Cylindrical, and Ellipsoidal Harmonics, with Applications to
   Problems in Mathematical Physics\/}
   (Ginn, Boston, 1893)
   [reprinted: (Dover, Mineola, N.Y., 2003)].
\bibitem{Hobs96}
   E.W.\ Hobson,
   Phil.\ Trans.\ R.\ Soc.\ Lond.\ A 187 (1896) 443.
\bibitem{Wang04}
   A.\ Wangerin,
   \emph{Theorie der Kugelfunktionen und der verwandten Funktionen,
   insbesondere der Lam{\'e}'schen und Bessel'schen},
   in \emph{Encyklop{\"a}die der mathematischen Wissenschaften\/},
   vol.\ 2.1 (Teubner, Leipzig, 1904), p.\ 695.
\bibitem{Barn07}
   E.W.\ Barnes,
   Quart.\ J.\ Pure Appl.\ Math.\ 39 (1907) 97.
   The associated Legendre function of the second kind defined in
   that work differs from the counterpart function of
   Hobson \cite{Hobs31} used in the present paper. The relationship 
   between the two functions is: 
   $[Q_{\nu}^{\mu}(z)]_{\mathrm{Barnes}}
   =\{\mathrm{e}^{-\mathrm{i}\pi\mu}\sin[\pi(\nu+\mu)]/\sin(\pi\nu)\}
   Q_{\nu}^{\mu}(z)$.
\bibitem{Wang21}
   A.\ Wangerin,
   \emph{Theorie des Potentials und der Kugelfunktionen\/}, vol.\ 2
   (de Gruyter, Berlin, 1921).
\bibitem{Hobs31}
   E.W.\ Hobson,
   \emph{The Theory of Spherical and Ellipsoidal Harmonics\/}
   (Cambridge University Press, Cambridge, 1931)
   [reprinted: (Chelsea, New York, 1955)].
\bibitem{Snow52}
   Ch.\ Snow,
   \emph{Hypergeometric and Legendre Functions with Applications to
   Integral Equations of Potential Theory\/}, 2nd ed.
   (National Bureau of Standards, Washington, D.C., 1952).  
   The associated Legendre functions defined in that book differ from
   the counterpart functions of Hobson \cite{Hobs31} used in the 
   present paper. The relationships between the two sets of functions 
   are: $[P_{\nu}^{\mu}(z)]_{\mathrm{Snow}}
   =[\Gamma(\nu+\mu+1)/\Gamma(\nu-\mu+1)]P_{\nu}^{-\mu}(z)$ and
   $[Q_{\nu}^{\mu}(z)]_{\mathrm{Snow}}
   =\mathrm{e}^{-\mathrm{i}\pi\mu}\cos(\pi\mu)Q_{\nu}^{\mu}(z)$.
\bibitem{Erde53}
   A.\ Erd{\'e}lyi (ed.),
   \emph{Higher Transcendental Functions\/}, vol.\ 1 
   (McGraw-Hill, New York, 1953), chap.\ III.
\bibitem{Mors53}
   P.M.\ Morse and H.\ Feshbach,
   \emph{Methods of Theoretical Physics\/}
   (McGraw-Hill, New York, 1953).
\bibitem{Lens54}
   J.\ Lense,
   \emph{Kugelfunktionen\/}, 2nd ed.\
   (Geest \& Portig, Leipzig, 1954).
\bibitem{Robi57}
   L.\ Robin,
   \emph{Fonctions Sph{\'e}riques de Legendre et Fonctions
   Sph{\'e}ro{\"{\i}}dales\/}, vol.\ 1
   (Gauthier-Villars, Paris, 1957).
\bibitem{Robi58}
   L.\ Robin,
   \emph{Fonctions Sph{\'e}riques de Legendre et Fonctions
   Sph{\'e}ro{\"{\i}}dales\/}, vol.\ 2
   (Gauthier-Villars, Paris, 1958).
\bibitem{Robi59}
   L.\ Robin,
   \emph{Fonctions Sph{\'e}riques de Legendre et Fonctions
   Sph{\'e}ro{\"{\i}}dales\/}, vol.\ 3
   (Gauthier-Villars, Paris, 1959).
\bibitem{Jahn60}
   E.\ Jahnke, F.\ Emde, and F.\ L{\"o}sch,
   \emph{Tafeln h{\"o}herer Funktionen\/}, 6th ed.\
   (Teubner, Stuttgart, 1960).
\bibitem{Krat60}
   A.\ Kratzer and W.\ Franz,
   \emph{Transzendente Funktionen\/}
   (Akademische Verlagsgesellschaft, Leipzig, 1960), chap.\ 5.
\bibitem{Steg65}
   I.A.\ Stegun,
   \emph{Legendre functions},
   in \emph{Handbook of Mathematical Functions\/},
   M.\ Abramowitz and I.A.\ Stegun (eds.)
   (Dover, New York, 1965), p.\ 331.
\bibitem{Magn66}
   W.\ Magnus, F.\ Oberhettinger, and R.P.\ Soni,
   \emph{Formulas and Theorems for the Special Functions of
   Mathematical Physics\/}, 3rd ed.\
   (Springer, Berlin, 1966).
\bibitem{MacR67}
   T.M.\ MacRobert,
   \emph{Spherical Harmonics\/}, 3rd ed.\
   (Pergamon, Oxford, 1967).
\bibitem{Grad94}
   I.S.\ Gradshteyn and I.M.\ Ryzhik,
   \emph{Table of Integrals, Series, and Products\/}, 5th ed.\
   (Academic, San Diego, 1994).
\bibitem{Temm96}
   N.M.\ Temme,
   \emph{Special Functions.\ An Introduction to the Classical
   Functions of Mathematical Physics\/}
   (Wiley, New York, 1996), chap.\ 8.
\bibitem{Prud83}
   A.P.\ Prudnikov, Yu.A.\ Brychkov, and O.I.\ Marichev,
   \emph{Integrals and Series. Special Functions\/}
   (Nauka, Moscow, 1983) (in Russian).
\bibitem{Prud03}
   A.P.\ Prudnikov, Yu.A.\ Brychkov, and O.I.\ Marichev,
   \emph{Integrals and Series. Special Functions. Supplementary
   Chapters\/}, 2nd ed.\
   (Fizmatlit, Moscow, 2003) (in Russian).
\bibitem{Bryc06}
   Yu.A.\ Brychkov,
   \emph{Special Functions. Derivatives, Integrals, Series, and Other
   Formulae\/}
   (Fizmatlit, Moscow, 2006) (in Russian).
\bibitem{Davi65}
   P.J.\ Davis,
   \emph{Gamma function and related functions},
   in \emph{Handbook of Mathematical Functions\/},
   M.\ Abramowitz and I.A.\ Stegun (eds.)
   (Dover, New York, 1965), p.\ 253.
\bibitem{Brow95}
   G.J.N.\ Brown,
   J.\ Phys.\ A 28 (1995) 2297.
\bibitem{Wats18}
   G.N.\ Watson,
   Proc.\ Lond.\ Math.\ Soc.\ 17 (1918) 241.
\end{thebibliography}
\end{document}